\documentclass[10pt]{amsart}
\usepackage{amsmath}
\usepackage{amsxtra}
\usepackage{amscd}
\usepackage{amsthm}
\usepackage{amsfonts}
\usepackage{amssymb}
\usepackage{eucal}

\newtheorem{cor}[subsection]{Corollary}
\newtheorem{lem}[subsection]{Lemma}
\newtheorem{prop}[subsection]{Proposition}

\newtheorem{conj}[subsection]{Conjecture}
\newtheorem{thm}[subsection]{Theorem}
\newtheorem{defn}[subsection]{Definition}

\theoremstyle{definition}

\theoremstyle{remark}

\newcommand{\thmref}[1]{Theorem~\ref{#1}}
\newcommand{\secref}[1]{Sect.~\ref{#1}}
\newcommand{\lemref}[1]{Lemma~\ref{#1}}
\newcommand{\propref}[1]{Proposition~\ref{#1}}
\newcommand{\corref}[1]{Corollary~\ref{#1}}

\newcommand{\nc}{\newcommand}
\nc{\renc}{\renewcommand}
\nc{\ssec}{\subsection}
\nc{\sssec}{\subsubsection}
\nc{\on}{\operatorname}

\nc\ol{\overline}
\nc\wt{\widetilde}
\nc\wh{\widehat}
\nc\tboxtimes{\wt{\boxtimes}}

\nc{\BA}{{\mathbb{A}}}
\nc{\BC}{{\mathbb{C}}}
\nc{\BG}{{\mathbb{G}}}
\nc{\BH}{{\mathbb{H}}}
\nc{\BM}{{\mathbb{M}}}
\nc{\BN}{{\mathbb{N}}}
\nc{\BK}{{\mathbb{K}}}
\nc{\BP}{{\mathbb{P}}}
\nc{\BR}{{\mathbb{R}}}
\nc{\BX}{{\mathbb{X}}}
\nc{\BV}{{\mathbb{V}}}
\nc{\BU}{{\mathbb{U}}}
\nc{\BW}{{\mathbb{W}}}
\nc{\BZ}{{\mathbb{Z}}}
\nc{\BS}{{\mathbb{S}}}

\nc{\CA}{{\mathcal{A}}}
\nc{\CB}{{\mathcal{B}}}

\nc{\CE}{{\mathcal{E}}}
\nc{\CF}{{\mathcal{F}}}
\nc{\CG}{{\mathcal{G}}}
\nc{\CL}{{\mathcal{L}}}
\nc{\CC}{{\mathcal{C}}}
\nc{\CM}{{\mathcal{M}}}
\nc{\CN}{{\mathcal{N}}}
\nc{\CO}{{\mathcal{O}}}
\nc{\CP}{{\mathcal{P}}}
\nc{\CQ}{{\mathcal{Q}}}
\nc{\CR}{{\mathcal{R}}}
\nc{\CS}{{\mathcal{S}}}
\nc{\CT}{{\mathcal{T}}}
\nc{\CU}{{\mathcal{U}}}
\nc{\CV}{{\mathcal{V}}}
\nc{\CK}{{\mathcal{K}}}
\nc{\CW}{{\mathcal{W}}}
\nc{\CZ}{{\mathcal{Z}}}

\nc{\cM}{{\check{\mathcal M}}{}}
\nc{\csM}{{\check{\mathcal A}}{}}
\nc{\oM}{{\overset{\circ}{\mathcal M}}{}}
\nc{\obM}{{\overset{\circ}{\mathbf M}}{}}
\nc{\oCA}{{\overset{\circ}{\mathcal A}}{}}
\nc{\obA}{{\overset{\circ}{\mathbf A}}{}}
\nc{\ooM}{{\overset{\circ}{M}}{}}
\nc{\osM}{{\overset{\circ}{\mathsf M}}{}}
\nc{\vM}{{\overset{\bullet}{\mathcal M}}{}}
\nc{\nM}{{\underset{\bullet}{\mathcal M}}{}}
\nc{\oD}{{\overset{\circ}{\mathcal D}}{}}
\nc{\obD}{{\overset{\circ}{\mathbf D}}{}}
\nc{\oA}{{\overset{\circ}{\mathbb A}}{}}
\nc{\op}{{\overset{\bullet}{\mathbf p}}{}}
\nc{\cp}{{\overset{\circ}{\mathbf p}}{}}
\nc{\oU}{{\overset{\bullet}{\mathcal U}}{}}
\nc{\oZ}{{\overset{\circ}{\mathcal Z}}{}}
\nc{\ofZ}{{\overset{\circ}{\mathfrak Z}}{}}
\nc{\oF}{{\overset{\circ}{\fF}}}

\nc{\fa}{{\mathfrak{a}}}
\nc{\fb}{{\mathfrak{b}}}
\nc{\fg}{{\mathfrak{g}}}
\nc{\fgl}{{\mathfrak{gl}}}
\nc{\fh}{{\mathfrak{h}}}
\nc{\fj}{{\mathfrak{j}}}
\nc{\fm}{{\mathfrak{m}}}
\nc{\fn}{{\mathfrak{n}}}
\nc{\fu}{{\mathfrak{u}}}
\nc{\fp}{{\mathfrak{p}}}
\nc{\fr}{{\mathfrak{r}}}
\nc{\fs}{{\mathfrak{s}}}
\nc{\fsl}{{\mathfrak{sl}}}
\nc{\hsl}{{\widehat{\mathfrak{sl}}}}
\nc{\hgl}{{\widehat{\mathfrak{gl}}}}
\nc{\hg}{{\widehat{\mathfrak{g}}}}
\nc{\chg}{{\widehat{\mathfrak{g}}}{}^\vee}
\nc{\hn}{{\widehat{\mathfrak{n}}}}
\nc{\chn}{{\widehat{\mathfrak{n}}}{}^\vee}

\nc{\fA}{{\mathfrak{A}}}
\nc{\fB}{{\mathfrak{B}}}
\nc{\fD}{{\mathfrak{D}}}
\nc{\fE}{{\mathfrak{E}}}
\nc{\fF}{{\mathfrak{F}}}
\nc{\fG}{{\mathfrak{G}}}
\nc{\fK}{{\mathfrak{K}}}
\nc{\fL}{{\mathfrak{L}}}
\nc{\fM}{{\mathfrak{M}}}
\nc{\fN}{{\mathfrak{N}}}
\nc{\fP}{{\mathfrak{P}}}
\nc{\fU}{{\mathfrak{U}}}
\nc{\fV}{{\mathfrak{V}}}
\nc{\fZ}{{\mathfrak{Z}}}

\nc{\bb}{{\mathbf{b}}}
\nc{\bc}{{\mathbf{c}}}
\nc{\bd}{{\mathbf{d}}}
\nc{\be}{{\mathbf{e}}}
\nc{\bj}{{\mathbf{j}}}
\nc{\bn}{{\mathbf{n}}}
\nc{\bp}{{\mathbf{p}}}
\nc{\bq}{{\mathbf{q}}}
\nc{\bF}{{\mathbf{F}}}
\nc{\bu}{{\mathbf{u}}}
\nc{\bv}{{\mathbf{v}}}
\nc{\bx}{{\mathbf{x}}}
\nc{\bs}{{\mathbf{s}}}
\nc{\by}{{\mathbf{y}}}
\nc{\bw}{{\mathbf{w}}}
\nc{\bA}{{\mathbf{A}}}
\nc{\bK}{{\mathbf{K}}}
\nc{\bI}{{\mathbf{I}}}
\nc{\bB}{{\mathbf{B}}}
\nc{\bG}{{\mathbf{G}}}
\nc{\bC}{{\mathbf{C}}}
\nc{\bD}{{\mathbf{D}}}
\nc{\bP}{{\mathbf{P}}}
\nc{\bH}{{\mathbf{H}}}
\nc{\bM}{{\mathbf{M}}}
\nc{\bN}{{\mathbf{N}}}
\nc{\bV}{{\mathbf{V}}}
\nc{\bU}{{\mathbf{U}}}
\nc{\bL}{{\mathbf{L}}}
\nc{\bT}{{\mathbf{T}}}
\nc{\bW}{{\mathbf{W}}}
\nc{\bX}{{\mathbf{X}}}
\nc{\bY}{{\mathbf{Y}}}
\nc{\bZ}{{\mathbf{Z}}}
\nc{\bS}{{\mathbf{S}}}

\nc{\sA}{{\mathsf{A}}}
\nc{\sB}{{\mathsf{B}}}
\nc{\sC}{{\mathsf{C}}}
\nc{\sD}{{\mathsf{D}}}
\nc{\sF}{{\mathsf{F}}}
\nc{\sG}{{\mathsf{G}}}
\nc{\sK}{{\mathsf{K}}}
\nc{\sM}{{\mathsf{M}}}
\nc{\sH}{{\mathsf{H}}}
\nc{\sO}{{\mathsf{O}}}
\nc{\sQ}{{\mathsf{Q}}}
\nc{\sP}{{\mathsf{P}}}
\nc{\sZ}{{\mathsf{Z}}}
\nc{\sfp}{{\mathsf{p}}}
\nc{\sr}{{\mathsf{r}}}
\nc{\sg}{{\mathsf{g}}}
\nc{\sff}{{\mathsf{f}}}
\nc{\sfb}{{\mathsf{b}}}
\nc{\sfc}{{\mathsf{c}}}
\nc{\sd}{{\mathsf{d}}}

\nc{\tA}{{\widetilde{\mathbf{A}}}}
\nc{\tB}{{\widetilde{\mathcal{B}}}}
\nc{\tg}{{\widetilde{\mathfrak{g}}}}
\nc{\tG}{{\widetilde{G}}}
\nc{\TM}{{\widetilde{\mathbb{M}}}{}}
\nc{\tO}{{\widetilde{\mathsf{O}}}{}}
\nc{\tU}{{\widetilde{\mathfrak{U}}}{}}
\nc{\TZ}{{\tilde{Z}}}
\nc{\tx}{{\tilde{x}}}
\nc{\tbv}{{\tilde{\bv}}}
\nc{\tfP}{{\widetilde{\mathfrak{P}}}{}}
\nc{\tz}{{\tilde{\zeta}}}
\nc{\tmu}{{\tilde{\mu}}}

\nc{\one}{{\mathbf{1}}}
\nc{\two}{{\mathbf{t}}}

\nc{\Rep}{{\mathop{\operatorname{\rm Rep}}}}
\nc{\Tot}{{\mathop{\operatorname{\rm Tot}}}}
\nc{\Ker}{{\mathop{\operatorname{\rm Ker}}}}
\nc{\Hilb}{{\mathop{\operatorname{\rm Hilb}}}}
\nc{\End}{{\mathop{\operatorname{\rm End}}}}
\nc{\Ext}{{\mathop{\operatorname{\rm Ext}}}}
\nc{\Hom}{\on{\Hom}}
\nc{\CHom}{{\mathop{\operatorname{{\mathcal{H}}\it om}}}}
\nc{\GL}{{\mathop{\operatorname{\rm GL}}}}
\nc{\gr}{{\mathop{\operatorname{\rm gr}}}}
\nc{\Id}{{\mathop{\operatorname{\rm Id}}}}
\nc{\de}{{\mathop{\operatorname{\rm def}}}}
\nc{\length}{{\mathop{\operatorname{\rm length}}}}
\nc{\supp}{{\mathop{\operatorname{\rm supp}}}}

\nc{\Cliff}{{\mathsf{Cliff}}}
\nc{\Gr}{{\on{Gr}}}
\nc{\bGr}{{\mathbf Gr}}
\nc{\bSt}{{\mathbf St}}
\nc{\bFl}{{\mathbf Fl}}
\nc{\Fl}{\on{Fl}}
\nc{\Fib}{{\mathsf{Fib}}}
\nc{\Coh}{{\mathsf{Coh}}}
\nc{\FCoh}{{\mathsf{FCoh}}}
\nc{\BVect}{{\mathbb Vect}}
\nc{\BSet}{{\mathbb Set}}
\nc{\BRep}{{\mathbb Rep}}
\nc{\bSet}{{\mathbf Set}}

\nc{\reg}{{\text{\rm reg}}}

\nc{\cplus}{{\mathbf{C}_+}}
\nc{\cminus}{{\mathbf{C}_-}}
\nc{\cthree}{{\mathbf{C}_*}}
\nc{\Qbar}{{\bar{Q}}}

\nc{\bh}{{\bar{h}}}
\nc{\bOmega}{{\overline{\Omega}}}

\nc{\seq}[1]{\stackrel{#1}{\sim}}

\title[Algebraic groups over
a two-dimensional local field]{Representations of algebraic 
groups over a 2-dimensional local field}

\author{David Kazhdan and Dennis Gaitsgory}

\address{\newline
D.K.: Einstein Institute of Mathematics, the Hebrew
University of Jerusalem, Givat Ram, Jerusalem, 91904, Israel; \newline
D.G.: Department of Mathematics, The University of Chicago,
5734 University Ave., Chicago, IL, 60637, USA.} 

\email{kazhdan@math.huji.ac.il; gaitsgde@math.uchicago.edu}

\begin{document}

\begin{abstract}

We introduce a categorical framework for the study
of representations of $G(\bF)$, where $G$ is a reductive group,
and $\bF$ is a 2-dimensional local field, i.e., $\bF=\bK((t))$,
where $\bK$ is a local field.

Our main result says that the space of functions on $G(\bF)$,
which is an object of a suitable category of representations
of $G(\bF)$ with the respect to the action of $G$ on itself
by left translations, becomes a representation of a certain 
central extension of $G(\bF)$, when we consider the action by
right translations.

\end{abstract}

\maketitle

\section*{Introduction}

\ssec{}

Let $\bK$ be a local field, and let us consider the field $\bF=\bK((t))$.
In his paper \cite{Kap}, Kapranov studied a certain representation 
of the group $G(\bF)$, where $G$ is a reductive group over $\bK$. 
He introduced a pro-vector space (we will denote it by $\BV$), 
on which the group $G(\bF)$ acts in a continuous way, and which 
may be thought of as an analogue of a principal series representation 
of usual $\fp$-adic groups.

Namely, $\BV$ is the (pro)-vector space of locally constant functions
with compact support on the set of $\bK$-points on the base affine
space of the loop group $G((t))$. (We remind that this base affine 
space is a principal $T$-bundle over the affine flag scheme 
corresponding to $G$, where $T$ is the Cartan subgroup.) 

Kapranov wrote down a certain algebra of endomorphisms of $\BV$ generated
by explicit intertwining operators, and proved that this algebra
is isomorphic to the (modified) double affine Hecke algebra. This
double affine Hecke algebra, which was introduced and studied by 
Cherednik, is clearly an object of great importance, and Kapranov's
work explained that it is related to groups over a 2-dimensional
field, such as $\bF$, in the same way as the usual affine Hecke 
algebra is related to $\fp$-adic groups.

\ssec{}

The present paper grew out of an attempt to put Kapranov's
ideas and results into a categorical framework. Our goal is to find
a category of smooth representations, let us denote it 
$\on{Rep}(\BG)$, which would 
contain Kapranov's representation (and its close relatives) as 
objects. Moreover, we want $\on{Rep}(\BG)$ to be abelian, so 
that the usual representation-theoretic questions, such as 
irreducibility, would make sense in it. We also want $\on{Rep}(\BG)$ 
to be as ``rigid'' or ``constrained'' as possible, and finally we want 
the definition of $\on{Rep}(\BG)$ to resemble the definition of
the category of smooth representations for usual $\fp$-adic groups.

After some categorical preliminaries in \secref{ind-pro}, we propose
a definition of $\on{Rep}(\BG)$ in \secref{categories of representations}.
A somewhat surprising feature of $\on{Rep}(\BG)$ is that, unlike 
most abelian categories that arise in representation theory,
the natural forgetful functor defined on $\on{Rep}(\BG)$ does
not map to the category of vector spaces, but rather to the
category $\BVect$  of pro-vector spaces. We remark that for the 
purposes of this paper, one could restrict to the subcategory
$\BVect^{\aleph_0}$ of projective systems indexed by countable
sets.

Let us recall that $\BVect$ is an abelian category, but it is not
semi-simple. In fact, the subcategory $\BVect^{\aleph_0}$ has
cohomological dimension $\leq 1$, and it can be visualized as follows:
An object of $\BVect^{\aleph_0}$ is called strict if it can be 
represented as a (filtering, countable) inverse system of vector 
spaces $\bV_i$, such that the arrows $\bV_i\to \bV_j$ are surjective. 
Strict objects of $\BVect^{\aleph_0}$ are the same as vector spaces 
endowed with a linear topology, with a countable fundamental 
system of neighbourhoods of zero, in which they are separated and 
complete. However, as is well-known, the 
category of such topological vector spaces is not abelian, which 
corresponds to the fact that strict objects of $\BVect^{\aleph_0}$ 
do not form an abelian subcategory. 

We do have a (left-exact) functor $limProj:\BVect\to Vect$, but the 
point of view taken in this paper, and which is largely borrowed from 
\cite{Kap}, is that we really have to work with the abelian category 
$\BVect$, and avoid taking projective limits.

\medskip

We justify the appearance of $\BVect$ by showing that $G(\bF)$ does
not have representations in any reasonable sense, unless we admit
pro-vector spaces.

\ssec{}

In \secref{induction} we show, generalizing the basic construction of
\cite{Kap}, how to produce non-trivial objects of $\on{Rep}(\BG)$.

Namely, let $\BH$ be a subgroup of $G(\bF)$, contained in the group 
$G[[t]](\bK)$ of 
$\bK$-points of the group $G[[t]]$ and equal to the preimage of a 
closed subgroup of
$(G[[t]]/G^i)(\bK)$ for some congruence subgroup $G^i$. 
Then, representations of $\BH$ on vector spaces,
as well as on pro-vector spaces, are notions that are easy to recover
from the usual representation theory of $\fp$-adic groups. 

We define two functors $\wt{i}^\BG_\BH$ and $i^\BG_\BH$ from $\on{Rep}(\BH,\BVect)$ 
to $\on{Rep}(\BG)$, such that the former is the right adjoint to
the tautological restriction functor $\on{Rep}(\BH,\BVect)\to \on{Rep}(\BG)$.
(In other words, $\wt{i}^\BG_\BH$ should be thought of as an ordinary induction
functor, whereas we think of $i^\BG_\BH$ as some sort of semi-infinite induction,
by analogy with the theory of modules over vertex algebras, cf. \cite{AG}.)

Kapranov's representation $\BV$ is exactly of the form $i^\BG_\BH(\BC)$,
where $\BH$ is the group of $\bK$-points of the unipotent radical of
the Iwahori subgroup of $G((t))$, and $\BC$ is the trivial representation.
In \secref{examples} we give a slight improvement of Kapranov's main result 
by showing that the (modified) Cherednik's algebra maps isomorphically
onto the ring $\on{End}_{\on{Rep}(\BG)}(\BV)$.

In addition, in \secref{examples} we discuss another series of examples of 
objects of $\on{Rep}(\BG)$ by applying the functor $i^\BG_\BH$ for
$\BH=G[[t]](\bK)$ and $G[[t]](\bK)$-representations, which are restrictions
of irreducible cuspidal representations of the $\fp$-adic group $G(\bK)$.
By analogy with the corresponding result in the theory of $\fp$-adic groups,
we conjecture that these objects are actually irreducible in $\on{Rep}(\BG)$,
and give some evidence in support of this conjecture.

\ssec{}

Finally, in \secref{Schwartz space} we formulate and prove the main
result of this paper. 

Suppose that the group $G$ acts on an algebraic variety $S$. In the theory
of $\fp$-adic groups one introduces the Schwartz space 
$\on{Funct}^{lc}_c(S(\bK))$ of locally constant compactly supported 
functions on the set of $\bK$-points of $S$, which is a smooth representation 
of the group $G(\bK)$.

The question that we want to address is whether one can define an analogue
of the Schwartz space, denoted in this paper by $M(\BS)$, which would be related 
to functions and/or distributions on the set of $\bF$-valued points of $S$.
Of course, one expects that $M(\BS)$ is an object of $\BVect$, underlying
a $G(\bF)$-representation.

It appears that the answer to this question is negative in the simplest example 
of $G=SL_2$ acting on the projective line, and the situation seems to be 
analogous 
to the problem of developing the theory of D-modules on loop spaces, 
cf. \cite{AG}. 

However, there are two important examples of $G$-varieties $S$, for which we
can define $M(\BS)$:

First, we consider the case of $S$ being the affine space $A^n$, with
the natural action of $GL_n$. We introduce a space $M(\BA^n)$ and show that
it is naturally an object in the category of representations of the group
$\wh{\BG L}_n$ (here $\wh{\BG L}_n$ is the group of $\bK$-points of the
canonical (i.e., Tate) central extension 
$1\to G_m\to \wh{GL}_n\to GL_n((t))\to 1$).

Next, we consider the case when the variety $S$ is isomorphic to the group
$G$ itself, with the action by left translations, and we construct an object
$M(\BG)\in \on{Rep}(\BG)$. Now the natural question to ask is, whether
the action of $G(\bF)$ on itself by {\it right translations} defines on
$M(\BG)$ another, commuting, structure of an object of $\on{Rep}(\BG)$.

The answer to this question is that the right action of $G(\bF)$
on $M(\BG)$ develops an anomaly (compare it with the main theorem
from \cite{AG}). Namely, $M(\BG)$ does carry a commuting action, but
of the group of $\bK$-points of the central extension 
$1\to G_m\to \wh{G}\to G((t))\to 1$ corresponding to the adjoint
action of $G$ on its Lie algebra.

\ssec{Acknowledgments}

We would like to thank S.~Arkhipov, I.~Cherednik,
P.~Etingof, V.Ginzburg, M.~Kapranov
for useful discussions and communications.   
We are grateful to the anonymous reviewer for valuable
comments on the previous version of the paper, and to
E.~Hrushovski for reading the revised version.

The research of D.G. is supported by the long-term fellowship at 
the Clay Mathematics Institute. He also wants to thank
the Mathematics Department of the Hebrew University of
Jerusalem, where the main part of this work was written.

\section{Preliminaries}  \label{ind-pro}

\ssec{}

We will work with inductive and projective limits of
objects of various categories. Let $I$ be a set. Recall that
$I$ is said to be filtering if it is endowed with a partial 
order, such that for any two elements $i_1,i_2\in I$, there
exists an element $i'\in I$ with $i'\geq i_1,i_2$. 

Let $I$ be a filtering set, which we can regard as a category, 
and $\Phi:i\mapsto S_i$ be a functor $I\to Set$. 
We will denote by $\underset{\longrightarrow}{lim}\, S_i$ the 
inductive limit of $\Phi$. 
In other words, $$\on{Hom}_{Set}(\underset{\longrightarrow}{lim}\, S_i,S)\simeq 
\on{Hom}_{\on{Functors}}(\Phi,\Phi_S),$$
where $\on{Functors}$ denotes the category of functors $I\to Set$,
and $\Phi_S$ is the ``constant'' functor corresponding to the set
$S$.

\medskip

Let $\CC$ be an arbitrary category. Recall from \cite{SGA4}
that the ind-completion of $\CC$, denoted $\on{Ind}(\CC)$, is the full
subcategory in the category of contravariant functors $\CC\to
Set$, which consists of objects (isomorphic to ones) of the form
$$X\mapsto \underset{\longrightarrow}{lim}\,\on{Hom}_{\CC}(X,X_i),$$
where $i\mapsto X_i$ is a functor $I\to \CC$ and $I$
is a filtering set;
we will denote by $"\underset{\longrightarrow}{lim}"\,X_i$ the
corresponding object of $\on{Ind}(\CC)$, which we will call 
``the direct limit of the system $X_i$''. By definition, 
$"\underset{\longrightarrow}{lim}"\,X_i(X)=\underset{\longrightarrow}{lim}\,\on{Hom}(X,X_i)$,
where the inductive limit is taken in the category of sets.

For example, let $Vect$ (resp., $Vect_0$) be the category of vector
spaces (resp., finite-dimensional vector spaces) over a given
ground field. We have $Vect\simeq \on{Ind}(Vect_0)$. 
(It is a good exercise to show $\on{Ind}(Vect)$
is NOT equivalent to $Vect$.)

For a cardinal $\aleph$, we will denote by $\on{Ind}^\aleph(\CC)$
the full subcategory of $\on{Ind}(\CC)$ obtained by imposing
the condition that the sets of indices that we are considering 
are of cardinality $\leq\aleph$.

\medskip

We have a canonical fully faithful embedding
$\CC\to \on{Ind}(\CC)$. The (partially defined) left adjoint
$\on{Ind}(\CC)\to \CC$, called the inductive limit, which we will
denote by $limInd$,  is always right-exact. 
We will say that $\CC$ is closed under inductive limits
(resp., inductive limits of cardinality $\leq\aleph$)
if the functor $limInd$ is defined 
on the entire $\on{Ind}(\CC)$
(resp., $\on{Ind}^\aleph(\CC)$). 
For example, it is easy to show that any category of the form $\on{Ind}(\CC)$,
(resp., $\on{Ind}^\aleph(\CC)$), where $\CC$ is another category, 
is always closed under inductive limits 
(resp., of cardinality $\leq\aleph$).

\medskip

\noindent {\bf Note on the terminology}:
Let us emphasize that  for a functor $I\to \CC:i\mapsto X_i$, we denote by 
$"\underset{\longrightarrow}{lim}"\,X_i$ the corresponding object of $\on{Ind}(\CC)$,
and call it the direct limit of the $X_i$'s, following \cite{SGA4}.
By contrast, if the object
$limInd("\underset{\longrightarrow}{lim}"\,X_i)\in \CC$ exists, we will call it
the inductive limit of the $X_i$'s, and denote it also by
$\underset{\longrightarrow}{lim}\,X_i$.

\medskip

The following simple assertion is useful:

\begin{lem} \label{ind vs representability}
Assume that $\CC$ is closed under inductive limits
of cardinality $\leq\aleph$, and $X\in \on{Ind}^\aleph(\CC)$. 
Then $X$ belongs to $\CC$ if and only if for every 
$"\underset{\longrightarrow}{lim}"X'_i=:X'\in \on{Ind}^\aleph(\CC)$,
the canonical arrow
$X\left(\underset{\longrightarrow}{lim}\, X'_i\right)\to 
\underset{\longleftarrow}{lim}\, \left(X(X'_i)\right)\simeq  
\on{Hom}_{\on{Ind}(\CC)}(X',X)$
is an isomorphism.
\end{lem}

\medskip

The pro-completion $\on{Pro}(\CC)$ (resp., $\on{Pro}^\aleph(\CC)$)
and the functor $limProj: \on{Pro}(\CC)\to \CC$
are defined in the same way by inverting the arrows, i.e.,
$\on{Pro}(\CC)=\left(\on{Ind}(\CC^o)\right)^o$, where the
superscript ``$o$'' means the opposite category.

\ssec{}

Suppose now that $\CC$ is an additive (resp., $\BC$-linear) category.
Then every object $F$
of $\on{Ind}(\CC)$, which is a priori a contravariant functor
$\CC\to Set$, lifts in a natural way to an additive functor
$\CC\to\{\text{Abelian groups}\}$ (resp., $\BC$-linear functor $\CC\to Vect$).

Indeed if for some $X_i\in \CC$,
$X="\underset{\longrightarrow}{lim}"\,X_i$, for the corresponding
$\on{Hom}$ sets we have:
$X(Y)=\underset{\longrightarrow}{lim}\on{Hom}(Y,X_i)$, and
this inductive limit of sets has a natural structure of an abelian group
(resp., $\BC$-vector space).

\medskip

Suppose now that $\CC$ is abelian. We will now give a simple 
criterion that establishes
ind-representability of functors in this case. Together with 
\lemref{ind vs representability} this provides a tool to prove representability
of various functors in the framework of abelian categories.

Assume that $\CC$ is such that for a given object the class if 
its subobjects is a set. Let $F:\CC\to \{\text{Abelian groups}\}$ 
be a contravariant left exact
functor. Suppose that there exists another functor $F'$ and a morphism of functors
$F\to F'$ such that $\forall X\in \CC$ the map $F(X)\to F'(X)$ is injective.

\begin{prop}  \label{prorep dom}
Assume that for $\CC$, $F$ and $F'$ as above,
the functor
$F'$ is ind-representable. Then the functor $F$ is also
ind-representable.
\end{prop}

\begin{proof} (Compare \cite{AM}, Corollary 2.8)

Let $Z="\underset{\longrightarrow}{lim}"\, Z_i$ be the object of
$\on{Ind}(\CC)$ ind-representing $F'$. For each index $i$ consider the functor
$F_i$ equal to $F\underset{F'}\times \on{Hom}(\cdot,Z_i)$. Each $F_i$
is also left exact and its map to $\on{Hom}(\cdot,Z_i)$ is an injection.
Obviously, $F(Y)=\underset{\longrightarrow}{lim}\, F_i(Y)$, so 
it enough to show that each $F_i$
is ind-representable. In other words, we can assume that $F'$ is representable
by an object $Z\in \CC$.

Consider the category of pairs
$(X\in \CC,\alpha:X\to Z)$, where $\alpha$ is an injective morphism in
$\CC$ such that the corresponding element in $F'(X)$ belongs to $F(X)$.
(Morphisms between $(X\in \CC,\alpha:X\to Z)$ and 
$(X'\in \CC,\alpha':X'\to Z)$ are maps $X\to X'$ in $\CC$, which
commute with the data of $\alpha$ and $\alpha'$.)
This category is obviously discrete, and it is small due to our assumption
on $\CC$. This resulting poset is filtering and is endowed with a 
functor to $\CC$, i.e., $(X\in \CC,\alpha:X\to Z)\mapsto X$. 

Let $W\in \on{Ind}(\CC)$ be the direct limit of this system. We claim that
$W$ ind-represents the functor $F$. Indeed, for $Y\in \CC$,
given an element in $\on{Hom}(Y,W)$ we have for some $X$ an element in $F(X)$
and a map $Y\to X$, which gives rise to an element of $F(Y)$.

And vice versa, given an element in $a_Y\in F(Y)$ consider the 
corresponding element 
$a'_Y\in F'(Y)$ and the resulting map $Y\to Z$. Let $X$ be the image of this map:
$$Y\twoheadrightarrow X\hookrightarrow Z.$$ It is enough to show that $a_Y$
belongs to the image of $F(X)$. Since $F$ is left exact, it is enough to show
that the image of $a_Y$ vanishes in $F\left(\on{ker}(Y\to X)\right)$. However, 
by assumption, 
the image of $a'_Y$ in $F'\left(\on{ker}(Y\to X)\right)$ is zero, which implies 
our assertion, since $F\to F'$ is injective.

\end{proof}

The following is also well-known (cf. \cite{AM}, Proposition 4.5):

\begin{lem}  \label{abelian}
If $\CC$ is abelian, then so is $\on{Ind}(\CC)$. The 
functor $limInd:\on{Ind}(\on{Ind}(\CC))\to 
\on{Ind}(\CC)$ is exact.
\end{lem}

Of course, assertions similar to the above ones hold when we replace $\on{Ind}$
by $\on{Pro}$.

\ssec{}   \label{set notation}

The following category will play an essential role in this paper: 
$$\BVect:=\on{Pro}(Vect)\simeq \on{Pro}(\on{Ind}(Vect_0)).$$
According to \lemref{abelian}, this is an abelian category. 

We will also consider the categories
$\bSet:=\on{Ind}(\on{Pro}(Set_0))$, and
$$\BSet:=\on{Ind}(\on{Pro}(\bSet))\simeq
\on{Ind}(\on{Pro}(\on{Ind}
(\on{Pro}(Set_0)))),$$ where $Set_0$
is the category of finite sets.

Note that the category 
$\on{Pro}(Set_0)$ is equivalent to the category of compact totally
disconnected topological spaces; let us denote this equivalence
by $\bY\mapsto \bY^{\on{top}}$. If $\bY="\underset{\longleftarrow}{lim}"\, Y_j$, 
then $\bY^{\on{top}}\simeq \bY=\underset{\longleftarrow}{lim}\, Y_j$, 
where the projective
limit is taken in the category of topological spaces.
For $\bX\in \bSet$ presented as a direct limit 
$"\underset{\longrightarrow}{lim}"\, \bX_i$ 
with  $\bX_i\in \on{Pro}(Set_0)$, set $\bX^{\on{top}}$ to be the topological 
space $\underset{\longrightarrow}{lim}\, \bX^{\on{top}}_i$ 
(where the inductive limit is again taken 
in the category of topological spaces). 

We will use the following terminology. We will call an object $\bX\in \bSet$ 
compact, if it belongs to $\on{Pro}(Set_0)$, and a morphism $\bX\to \bY$ 
in $\bSet$ proper if every base change by a compact object is compact.

We will call an object $\bX\in \bSet$ locally compact,
if it can be represented as a direct limit
$\bX="\underset{\longrightarrow}{lim}"\, \bX_i$, 
$\bX_i\in \on{Pro}(Set_0)$, where the maps $\bX_i\to \bX_j$ are 
such that the corresponding maps of topological spaces
$\bX^{\on{top}}_i\to \bX^{\on{top}}_j$ are open embeddings.

The full subcategory of $\bSet$ consisting of locally compact objects
is equivalent to the category $Top^{Hlctd}$ of Hausdorff locally compact
totally disconnected topological spaces. All objects of $\bSet$ that
are relevant for the purposes of this paper will be locally compact. 
Therefore, the reader may safely replace $\bSet$ by $Top^{Hlctd}$ and
$\BSet$ by $\on{Ind}(\on{Pro}(Top^{Hlctd}))$.

\medskip

Similarly, we will call an object $\BX\in \BSet$ bounded if it
actually belongs to $\on{Pro}(\bSet)$.

\ssec{}

Let $\CA$ be a monoidal category, i.e., we have a 
functor $\otimes:\CA\times \CA\to \CA$, a unit object $\one_\CA\in \CA$
and functorial isomorphisms
$$X\otimes (Y\otimes Z)\simeq (X\otimes Y)\otimes Z;\,\,\,
X\otimes \one_\CA\simeq X\simeq \one_\CA\otimes X,$$
obeying the usual axioms.
Note that in this case the categories $\on{Ind}(\CA)$ and $\on{Pro}(\CA)$ 
also possess natural monoidal structures.

If $\CC$ is another category, there is a standard notion of action of
$\CA$ on $\CC$, in which case we say that $\CC$ is a module category
over $\CA$. Namely, a module structure is a functor 
$\otimes: \CA\times \CC\to \CC$, and for $X,Y\in \CA$ and $V\in \CC$
functorial isomorphisms
$$(X\otimes Y)\otimes V\to X\otimes (Y\otimes V);\,\,\,
\one_\CA\otimes V\simeq V,$$
satisfying the natural axioms. In particular, for $X\in \CA$,
$V,W\in \CC$ we have a well-defined Hom set $\on{Hom}(X\otimes V,W)$.

By definition, a {\it pseudo-action} of $\CA$ on $\CC$ 
(or a structure on $\CC$ of a {\it pseudo-module} over $\CA$)
is a functor $\CA^o\times \CC^o\times \CC\to Set$, denoted $\CHom(\cdot\otimes
\cdot,\cdot)$, and a morphism of functors: for $X,Y\in \CA$, $V,U,W\in \CC$
$$\CHom(X\otimes V,W)\times \CHom(Y\otimes U,V)\Rightarrow
\CHom((X\otimes Y)\otimes U,W),$$
and a functorial isomorphism $\CHom(\one_\CA\otimes V,W)\simeq 
\on{Hom}_\CC(V,W)$, such that the following compatibility conditions hold: 

\smallskip

\noindent For $X,Y,Z\in \CA$, $V,U,W,Q\in \CC$, the arrows
\begin{align*}
&\CHom(X\otimes W,Q)\times \CHom(Y\otimes V,W) \times \CHom(Z\otimes U,V)
\to \\
&\to \CHom(X\otimes W,Q)\times \CHom((Y\otimes Z)\otimes U,W)\to
\CHom((X\otimes (Y\otimes Z))\otimes U,Q) \,\text{  and} \\
&\CHom(X\otimes W,Q)\times \CHom(Y\otimes V,W) 
\times \CHom(Z\otimes U,V)\to \\
&\to \CHom((X\otimes Y)\otimes V,Q)\otimes \CHom(Z\otimes U,V) \to
\CHom(((X\otimes Y)\otimes Z)\otimes U,Q)
\end{align*}
coincide under the associativity isomorphism
$(X\otimes Y)\otimes Z\simeq X\otimes (Y\otimes Z)$, 
and for $U,V,W\in \CC$ and $X\in \CA$, the squares
$$
\CD
\CHom(X\otimes V,U)\times \on{Hom}_\CC(W,V) @>>>
\CHom(X\otimes V,U)\times \CHom(\one_\CA\otimes W,V) \\
@VVV   @VVV  \\
\CHom(X\otimes W,U) @>>>
\CHom((X\otimes\one_\CA)\otimes W,U)
\endCD
$$
and
$$
\CD
\on{Hom}_\CC(V,W)\times \CHom(X\otimes U,V) @>>>
\CHom(\one_\CA\otimes V,W)\times \CHom(X\otimes U,V) \\
@VVV    @VVV \\
\CHom(X\otimes U,W)  @>>> \CHom((\one_\CA\otimes X)\otimes U,W)
\endCD
$$
are commutative.

\medskip

Note that if $\CC$ is a pseudo-module over $\CA$, then $\CC^o$ is a 
pseudo-module over $\CA^{op}$, where the latter is the category $\CA$ with 
the opposite monoidal structure:
$${\mathcal Hom}(X\otimes V^o,W^o):={\mathcal Hom}(X\otimes W,V).$$

When $\CC$ is additive (resp., $\BC$-linear), we will rather
use the variant of the above definition, when we require that the
sets $\CHom(X\otimes U,V)$ have a structure of an
abelian group (resp., $\BC$-vector space), such that the 
natural transformations 
$\CHom(X\otimes U,V)\times \CHom(Y\otimes V,W)\Rightarrow
\CHom((X\otimes Y)\otimes U,W)$ and $\CHom(\one_\CA\otimes V,W)\simeq 
\on{Hom}_\CC(V,W)$ are bilinear (resp., linear).

For example, the category $Set$ is a monoidal 
via $X\otimes Y:=X\times Y$, and any category $\CC$ has a 
pseudo-module structure over $Set$ via $\CHom(X\otimes U,V):=Hom(U,V)^X$
for $X\in Set$, $U,V\in \CC$.

Of course, when $\CC$ is a module category over $\CA$, it
acquires a pseudo-module structure by setting
$$\CHom(X\otimes U,V):=\on{Hom}_\CC(X\otimes U,V).$$

In what follows we will say that an element $\phi\in 
\CHom(X\otimes U,V)$ defines an action $X\times U\to V$.

\ssec{}

Let us now analyze how pseudo-actions behave when we Ind- and Pro-
complete our categories.

First, we claim that if $\CA$ pseudo-acts on $\CC$, then so do
$\on{Ind}(\CA)$ and $\on{Pro}(\CA)$. Indeed, if
$X\in \on{Ind}(\CA)$ (resp., $X\in \on{Pro}(\CA)$) is
$"\underset{\longrightarrow}{lim}"\,X_i$ 
(resp., $"\underset{\longleftarrow}{lim}"\,X_i$),
we set $\CHom(X\otimes V,W)=\underset{\longleftarrow}{lim} \,
\CHom(X_i\otimes V,W)$
(resp., $\CHom(X\otimes V,W)=\underset{\longrightarrow}{lim} \,
\CHom(X_i\otimes V,W)$).
It is easy to see that this definition is independent of the way we
represent $X$ as a direct (resp., inverse) limit.

Also, if $\CC$ has a pseudo-module structure over $\CA$, so do
$\on{Ind}(\CC)$ and $\on{Pro}(\CC)$. Indeed, for $V,W\in \on{Ind}(\CC)$
equal to $"\underset{\longrightarrow}{lim}"\,V_i$ and
$"\underset{\longrightarrow}{lim}"\,W_j$, 
(resp., $"\underset{\longleftarrow}{lim}"\,V_i$ and 
$"\underset{\longleftarrow}{lim}"\,W_j$), we set $\CHom(X\otimes V,W)$ to be
$$(\underset{i}{\underset{\longleftarrow}{lim}}) 
(\underset{j}{\underset{\longrightarrow}{lim}})\,\,\CHom(X\otimes V_i,W_j)
\text{  and  }
(\underset{j}{\underset{\longleftarrow}{lim}}) 
(\underset{i}{\underset{\longrightarrow}{lim}})\,\,\CHom(X\otimes V_i,W_j),$$
respectively. One can easily see that this definition is independent of the
presentation of $V$ and $W$ as directs (resp., inverse) limits. 

\medskip

Now, we obtain that there are two pseudo-actions of $\on{Ind}(\CA)$
on $\on{Ind}(\CC)$. One is (which we will call "naive")
when we first consider the pseudo-action of $\on{Ind}(\CA)$ on $\CC$
and then produce from it the corresponding pseudo-action on $\on{Ind}(\CC)$.
The other is when we first consider
the pseudo-action of $\CA$ on $\on{Ind}(\CC)$ and then produce from it
the corresponding pseudo-action of $\on{Ind}(\CA)$. Unless specified
otherwise, in the sequel we will use the pseudo-action
of the second kind. Note that we have a canonical map
$\CHom(X\otimes V,W)_{naive}\to \CHom(X\otimes V,W)$.
In concrete terms, if $X="\underset{\longrightarrow}{lim}"\,X_k$,
$"\underset{\longrightarrow}{lim}"\,V_i$ and
$"\underset{\longrightarrow}{lim}"\,W_j$, we have:
\begin{align*}
&\CHom(X\otimes V,W)_{naive}=\underset{i}{(\underset{\longleftarrow}{lim})}
\underset{j}{(\underset{\longrightarrow}{lim})}
\underset{k}{(\underset{\longleftarrow}{lim})}\,\,
\CHom(X_k\otimes V_i,W_j); \\
&\CHom(X\otimes V,W)=\underset{k}{(\underset{\longleftarrow}{lim})}
\underset{i}{(\underset{\longleftarrow}{lim})}
\underset{j}{(\underset{\longrightarrow}{lim})}\,\,
\CHom(X_k\otimes V_i,W_j).
\end{align*}

\medskip

For example, by taking $\CC=\CA$, the canonical action of
$\on{Ind}(\CA)$ on itself corresponding to the monoidal structure
coincides with the pseudo-action described above coming from
the action on $\CA$ on itself.

Similarly, we obtain the corresponding notions
concerning the pseudo-action of $\on{Ind}(\CA)$ on $\on{Pro}(\CC)$.

\medskip

The situation with the pseudo-actions of $\on{Pro}(\CA)$ is the opposite.
The naive pseudo-module structure on $\on{Ind}(\CC)$ is obtained when
we first consider the pseudo-action of $\CA$ on $\on{Ind}(\CC)$,
and then produce from it a pseudo-action of $\on{Pro}(\CA)$.
The pseudo-module structure that we will normally consider is
is obtained by first considering the pseudo-action of $\on{Pro}(\CA)$
on $\CC$, and then producing from it the corresponding
pseudo-action on $\on{Ind}(\CC)$. As before, we have a canonical map
$\CHom(X\otimes V,W)_{naive}\to \CHom(X\otimes V,W)$, and for
$V="\underset{\longrightarrow}{lim}"\,V_i$, 
$W="\underset{\longrightarrow}{lim}"\,W_j$
and $X="\underset{\longleftarrow}{lim}"\,X_k$
\begin{align*}
&\CHom(X\otimes V,W)_{naive}=\underset{k}{(\underset{\longrightarrow}{lim})}
\underset{i}{(\underset{\longleftarrow}{lim})}
\underset{j}{(\underset{\longrightarrow}{lim})}\,\,
\CHom(X_k\otimes V_i,W_j); \\
&\CHom(X\otimes V,W)=\underset{i}{(\underset{\longleftarrow}{lim})}\underset{j}
{(\underset{\longrightarrow}{lim})}
\underset{k}{(\underset{\longrightarrow}{lim})}\,\,
\CHom(X_k\otimes V_i,W_j).
\end{align*}

In a similar way, we obtain the two pseudo-actions of
$\on{Pro}(\CA)$ on $\on{Pro}(\CC)$.
As above, for $\CC=\CA$ this canonical pseudo-action coincides
with the action corresponding to the monoidal structure on $\on{Pro}(\CA)$.

\medskip

Finally, we see that there are 3 possible pseudo-actions 
of $\on{Ind}\on{Pro}(\CA)$
on $\on{Ind}(\CC)$. The one that we will consider is "the biggest":
we will first consider the pseudo-action of $\on{Pro}(\CA)$ on
$\CC$, then produce from it the pseudo-action of $\on{Pro}(\CA)$
on $\on{Ind}(\CC)$, and then the pseudo-action of $\on{Ind}\on{Pro}(\CA)$ on
$\on{Ind}(\CC)$.

Explicitly, if $X\in \on{Ind}\on{Pro}(\CA)$ is 
$\underset{k}{"\underset{\longrightarrow}{lim}"}
(\underset{l}{"\underset{\longleftarrow}{lim}"\, X^k_l})$,
$V="\underset{\longrightarrow}{lim}"\,V_i$ and 
$W="\underset{\longrightarrow}{lim}"\,W_j$, then
$$\CHom(X\otimes V,W)=
\underset{k}{\underset{\longleftarrow}{lim}}\, 
\underset{i}{\underset{\longleftarrow}{lim}}\, 
\underset{l}{\underset{\longrightarrow}{lim}}\,
\underset{j}{\underset{\longrightarrow}{lim}}\, \CHom(X^k_l\otimes V_i,W_j).$$

By inverting the arrows in $\CC$ we obtain the corresponding pseudo-action of
$\on{Ind}\on{Pro}(\CA)$ on $\on{Pro}(\CC)$.

\ssec{}   \label{ex}

Let us consider our main examples. Let $\CA=Set_0$,
and $\CC=Vect_0$. Then for
$\bX\in \bSet=
\on{Ind}(\on{Pro}(Set_0))$, 
$\bV,\bW\in Vect=\on{Ind}(Vect_0)$,
we obtain the notion of an action $\bX\times\bV\to \bW$.
However, it is easy to see that such an action 
is the same as a continuous map $\bX^{\on{top}}\times \bV\to \bW$, 
linear in $\bV$ and $\bW$, where $\bV$ and $\bW$ are endowed 
with the discrete topology, and $\bX^{\on{top}}$ is as in \secref{set notation}.

\medskip

Now set $\CA=\bSet=\on{Ind}(\on{Pro}(Set_0))$ and $\CC=Vect$.
We obtain a pseudo-module
structure on $\BVect$ with respect to $\BSet$.

Let us write down the last notion in more concrete terms.
First, let $\BX$ be an object of
$\on{Pro}(\bSet)$,
and $\BV,\BW$ be two objects of $\BVect$.
An action $\phi:\BX\times \BV\to \BW$ is the following data. Let
$\BX="\underset{\longleftarrow}{lim}"\,\bX_j$,
$\BV="\underset{\longleftarrow}{lim}"\,\bV_i$,
$\BW="\underset{\longleftarrow}{lim}"\,\bW_{i'}$, with
$\bX_j\in \bSet$, $\bV_i,\bW_{i'}\in Vect$.
Then for every $i'$ there must exist $i_0$, $j_0$
and a compatible system of action maps
$\phi_{j,i,i'}:\bX_j\times \bV_i\to \bW_{i'}$ defined for 
$i\geq i_0$, $j\geq j_0$.
Another compatibility condition is imposed: for $i'_1\geq i'_2$ the
corresponding diagrams
$$
\CD
\bX_{j_1}\times \bV_{i_1} @>{\phi_{j_1,i_1,i'_1}}>> \bW_{i'_1}  \\
@VVV    @VVV  \\
\bX_{j_2}\times \bV_{i_2} @>{\phi_{j_2,i_2,i'_2}}>> \bW_{i'_2}
\endCD
$$
must commute for $i_1$ and $j_1$ large enough.
Two action maps $\phi$ and $\psi$ coincide if for every
$i'$ the corresponding maps $\phi_{j,i,i'}$ and $\psi_{j,i,i'}$
coincide for $i$ and $j$ large enough.

If now $\BX$ is an object of
$\BSet$, equal to
$"\underset{\longrightarrow}{lim}"\,\BX^j$ and
$\BV,\BW\in \BVect$, an action
$\phi:\BX\times \BV\to \BW$ is a compatible system of
actions $\phi^j:\BX^j\times \BV\to \BW$.

\ssec{}   \label{weakly strict}

The following definition will be needed in the sequel. First, note that
we have an obvious functor from the category of sets (denoted $Set$) to
$\bSet$ via 
$$Set\simeq \on{Ind}(Set_0)\to \on{Ind}(\on{Pro}(Set_0))\simeq \bSet.$$
Let $\bX_1\to \bX_2$ be a map of objects of $\bSet$. We will say that
it is weakly surjective if for any $Y\in Set$, the map
$$\on{Hom}_{\bSet}(\bX_2,Y)\to \on{Hom}_{\bSet}(\bX_1,Y)$$
is injective. 

Note that if $\bX_1,\bX_2$ are locally compact, the above notion 
that a morphism $\bX_1\to \bX_2$ is weakly surjective is equivalent to
the condition that the corresponding map $\bX_1^{\on{top}}\to
\bX_2^{\on{top}}$ has dense image.

\begin{lem}
A map $\bX_1\to \bX_2$ in $\bSet$ is weakly surjective if
and only if for any $\bV,\bW\in Vect$, the map
$\CHom(\bX_2\otimes \bV,\bW)\to \CHom(\bX_1\otimes \bV,\bW)$ is injective.
\end{lem}

We will call an object $\BX\in \on{Pro}(\bSet)$ weakly strict if it
can be represented as $"\underset{\longleftarrow}{lim}"\,\bX_i$, where the
maps $\bX_j\to \bX_i$ are weakly surjective. 

\medskip

Note that if $\BX$ is weakly strict and $\bV,\bW\in Vect$, for any 
element $\phi\in \CHom(\BX\otimes \bV,\bW)$ we have well-defined
kernel and image of $\phi$. By definition, 
$\on{ker}(\phi)\subset \bV$ (resp., $\on{Im}(\phi)\subset\bW$) is
the maximal (resp., minimal) subspace $\bV'$ of $\bV$
(resp., $\bW'$ of $\bW$) having the property
that $\phi$ factors through an element 
$\phi'\in \CHom(\BX\otimes \bV/\bV',\bW')$ (resp., 
$\phi'\in \CHom(\BX\otimes \bV,\bW')$).

Indeed, both $\on{ker}(\phi)$ and $\on{Im}(\phi)$ are clearly 
well-defined when $\bX\in \bSet$. 
If now $\BX="\underset{\longleftarrow}{lim}"\,\bX_i$, with
weakly surjective maps, and $\phi$  comes from an element 
$\phi_i\in \CHom(\bX_i\otimes \bV,\bW)$,
then it is easy to see that $\on{ker}(\phi_i)\subset \bV$ and
$\on{Im}(\phi_i)\subset \bW$ are the sought-for subspaces.

\section{Categories of representations}   \label{categories of representations}

\ssec{}

In the abstract set-up of the previous section,
let us recall that an object $X\in \CA$ is called a monoid
(in the sense of the monoidal structure on $\CA$)
if we are given a (multiplication) map $X\otimes X\to X$ 
and a (unit) map $\one_\CA\to X$, which satisfy the
usual associativity and unit axioms.

In our examples, the monoidal structure on $\CA$
will be such that $X\otimes Y$ is isomorphic to 
the categorical direct product $X\times Y$.
Moreover, $\on{Hom}_{\CC}(X,\one_\CA)$ will be
a one-element set $\forall X\in \CC$.
Note that this property is inherited by both
$\on{Ind}(\CA)$ and $\on{Pro}(\CA)$. 

\medskip

In this case, it makes sense to speak about group-objects
in $\CA$: a monoid $X$ is called a group
if there exists a map $\gamma:X\to X$ (automatically unique)
such that the two compositions
$$X\overset{\Delta}\to X\times X\overset{\on{id}\times \gamma}
\longrightarrow X\times X\overset{\on{mult}}\longrightarrow X \text{ and }$$
$$X\overset{\Delta}\to X\times X\overset{\gamma\times \on{id}}
\longrightarrow X\times X\overset{\on{mult}}\longrightarrow X$$
are both equal to $X\to \one_\CA\to X$.

In the sequel we will only consider monoids, which are groups.

\ssec{}

If $\CC$ is another category with a pseudo-action of
$\CA$ and $X\in \CA$ is a monoid, a representation
of $X$ in $\CC$ is a pair $\Pi=(V,\rho)$, where
$V\in \CC$ and $\rho\in \CHom(X\otimes V,V)$, such that
the following two conditions hold:

\noindent Associativity: The image of $\rho\times\rho$ under
the associativity constraint
$$\CHom(X\otimes V,V)\otimes \CHom(X\otimes V,V)\to
\CHom((X\otimes X)\otimes V,V)$$
equals the image of $\rho$ under the map
$\CHom(X\otimes V,V)\to \CHom(X\otimes X)\otimes V,V)$ 
given by the multiplication $X\otimes X\to X$.

\noindent Unit: The image of $\rho$ in
$\CHom(\one_\CA\otimes V,V)$ under $\one_\CA\to X$
equals the identity element in
$\CHom(\one_\CA\otimes V,V)\simeq \on{Hom}(V,V)$.

\medskip

Representations of $X$ in $\CC$ form a category, which 
we will denote by $\on{Rep}(X,\CC)$. When $\CC$ is additive 
(resp., $\BC$-linear), the category $\on{Rep}(X,\CC)$ 
is additive (resp., $\BC$-linear) as well.

Assume now that 
$\CC$ is abelian and that for a fixed $X\in \CA$, the
functor $\CC^o\times \CC\to Set$ given by
$V,W\mapsto {\mathcal Hom}(X\otimes V,W)$ is left-exact
in both arguments.
\begin{lem}
Under the above circumstances the category
$\on{Rep}(X,\CC)$ is abelian and the natural
forgetful functor $\on{Rep}(X,\CC)\to \CC$ is exact.
\end{lem}
If $\CA$ is a monoidal category with a pseudo-action
on an abelian category $\CC$, such that the above
left-exactness condition is satisfied, then the same
holds for $\on{Ind}(\CA)$ (resp., $\on{Pro}(\CA)$)
pseudo-acting on $\on{Ind}(\CC)$ (resp., $\on{Pro}(\CC)$),
due to the fact that the functor $limInd$ (resp., $limProj$)
is exact (resp., left-exact) on the category of abelian
groups.

In particular, we obtain that this condition is satisfied
in our examples of $\CA=\bSet$, $\CC=Vect$ and $\CA=\BSet$,
$\CC=\BVect$.

\ssec{}  \label{condition star}

Set first $\CA=\bSet$, and $\CC=Vect$. Thus, for a group-object
$\bH\in \bSet$ the category $\on{Rep}(\bH,Vect)$
is the usual category of representations
of $\bH$ appearing in the theory of $\fp$-adic groups.
In other words, if $\bH$ is locally compact (cf. \secref{set notation})
and $\bH^{\on{top}}$ is the corresponding
topological group, then an object of $\on{Rep}(\bH,Vect)$ is 
the same as a smooth representation of $\bH^{\on{top}}$.

\medskip

If $\BH$ is a group-object of $\on{Pro}(\bSet)$,
we can consider its representations on $Vect$ and $\BVect$,
and the resulting categories will be denoted by 
$\on{Rep}(\BH,Vect)$ and $\on{Rep}(\BH,\BVect)$, respectively.

\medskip

We will say that $\BH\in \on{Pro}(\bSet)$ satisfies condition ($*$) if
it is weakly strict as an object of $\on{Pro}(\bSet)$, cf. \secref{weakly strict}.

The following assertion will play an important role in the sequel: 
\footnote{We would
like to thank E.~Hrushovski for pointing out the mistake in the previous
version of the paper, where \propref{representations of pro-groups},
was stated without the ($*$) assumption on $H$; in fact, he 
constructed a counterexample.}

\begin{prop}  \label{representations of pro-groups}
For $\BH$ satisfying ($*$), the categories $\on{Rep}(\BH,\BVect)$
and $\on{Pro}(\on{Rep}(\BH,Vect))$ are naturally equivalent.
\end{prop}

Note that the proof given below is valid when $\BH$ is a just a monoid
(not necessarily a group), satisfying condition $(*)$.

\begin{proof}

The functor in one direction:
$\sF:\on{Pro}\on{Rep}(\BH,Vect))\to \on{Rep}(\BH,\BVect)$ is evident;
moreover, it is easy to see that it is fully faithful.
Let us show that it admits a left adjoint. 

For $(\BV,\rho)\in\on{Rep}(\BH,\BVect)$,
let us write $\BV="\underset{\longleftarrow}{lim}"\, \bV_i$, where
the index $i$ runs over some filtering set $I$. Consider the category
of quadruples $(\bV',\rho',i,\alpha:\bV_i\to \bV')$, where
$(\bV',\rho')\in\on{Rep}(\BH,Vect)$, $i\in I$ and $\alpha$ is a map,
such that its image generates $\bV'$ as an $\BH$-representation and
the composition $\BV\to \bV_i\to \bV'$ is compatible with the $\BH$-actions.
A morphism between $(\bV_1,\rho_1',i_1,\alpha_1:\bV_{i_1}\to \bV'_1)$
and $(\bV_2,\rho_2',i_2,\alpha_2:\bV_{i_2}\to \bV'_2)$ is by definition
a relation $i_2\geq i_1$ and a map of $\BH$-representations $\bV'_1\to\bV'_2$,
such that the square
$$
\CD
\bV_{i_1} @>>> \bV_{i_2} \\
@V{\alpha_1}VV @V{\alpha_2}VV \\
\bV'_1 @>>> \bV'_2
\endCD
$$
commutes. The resulting category is evidently discrete, filtering and small.
By definition, we have a forgetful functor from this category to 
$\on{Rep}(\BH,Vect)$
that sends a quadruple $(\bV',\rho',i,\alpha)$ to $(\bV',\rho')$. Let us denote
by ${\mathsf G}(\BV,\rho)\in \on{Pro}(\on{Rep}(\BH,Vect))$ the resulting inverse 
limit. It is easy too see that the assignment $(\BV,\rho)\mapsto \sG(\BV,\rho)$
defines a functor left adjoint to $\sF$.

The fact that $\sF$ was fully-faifull means that the composition
$\sG\circ \sF$ is isomorphic to the identity functor. Thus, it remains
to see that for $(\BV,\rho)\in\on{Rep}(\BH,\BVect)$, the adjunction map
$(\BV,\rho)\to \sF\circ\sG(\BV,\rho)$ is an isomorphism. For that, it
suffices to show that if $\BV="\underset{\longleftarrow}{lim}"\, \bV_i$,
then for every $i$ there exists a vector space $\bV'_i$
underlying an object $\Pi=(\bV'_i,\rho_i)\in \on{Rep}(\BH,Vect)$,
such that the map $\BV\to \bV_i$ factors
as $\BV\to \bV'_i\to \bV_i$,
with the first arrow preserving the $\BH$-action. Indeed, this would show
that the map $(\BV,\rho)\to \sF\circ\sG(\BV,\rho)$ is always injective,
and combined with the fact that $\sG$ is right-exact, this implies 
that this map is an isomorphism.

Let $j$ be an index such that the map
$\BH\times \BV\overset{\on{act}}{\longrightarrow} \BV\to \bV_i$
factors as $\BH\times \BV\to \BH\times \bV_j\to \bV_i$.
Let us denote by $p_{j,i}$ the projection $\bV_j\to \bV_i$ and by
$\on{act}_{j,i}$ the map $\BH\times \bV_j\to \bV_i$.

By the definition of the action,
there exists another index $k$ such that the map
$\BH\times \BV\overset{\on{act}}{\longrightarrow} \BV\to \bV_j$
factors as
$$\BH\times \BV\to
\BH\times \bV_k\overset{\on{act}_{k,j}}\longrightarrow \bV_j,$$
and such that the diagram
$$
\CD
\BH\times \BH\times \bV_k @>{\on{mult}\times p_{k,j}}>> \BH\times
\bV_j \\
@V{\on{id}\times \on{act}_{k,j}}VV  @V{\on{act}_{j,i}}VV \\
\BH\times \bV_j @>{\on{act}_{j,i}}>>  \bV_i
\endCD
$$
is commutative, where $p_{k,j}$ denotes the projection
$\bV_k\to \bV_j$. 

Let $\bW'\subset \bV_j$ be the kernel of the map $\BH\times \bV_j
\overset{\on{act}_{j,i}}\longrightarrow \bV_i$, and let
$\bW''\subset \bV_j$ be the image of 
$\BH\times \bV_k\overset{\on{act}_{k,j}}\longrightarrow \bV_j$. 
The above kernel and image are well-defined due to the ($*$) assumption on 
$\BH$, cf. \secref{weakly strict}.

Set $\bV'_i$ to be the image of $\bW''$ in $\bV_j/\bW'$,
and let $\on{act}'$ denote the map
$\BH\times \bV_k\to \bV'_i$. We claim that there exists a 
unique map $\BH\times \bV'_i\to \bV'_i$, which makes the diagram
$$
\CD
\BH\times \BH\times \bV_k @>{\on{id}\times \on{act}'}>> \BH\times \bV'_i \\
@V{\on{mult}\times \on{id}}VV  @VVV   \\
\BH\times \bV_k @>{\on{act}'}>> \bV'_i 
\endCD
$$
commute. The commutativity of the diagram implies that the action
$\BH\times \bV'_i\to \bV'_i$ is unital and associative.

\medskip

To construct the sought-for map $\BH\times \bV'_i\to \bV'_i$,
let us write $\BH="\underset{\longleftarrow}{lim}"\, \bX_n$ with
weakly surjective maps. Let $n_0$ be an index such that the 
maps $\on{act}_{j,i}$ and
$\on{act}_{j,k}$ are defined on the level of $\bX_{n_0}$ (we will denote them
$\on{act}^{n_0}_{j,i}$ and $\on{act}^{n_0}_{k,j}$, respectively).

Let $n_1\geq n_0$ be an index such that the 
multiplication on $\BH$ gives rise to a map
$\on{mult}^{n_0}_{n_1,n_1}:\bX_{n_1}\times \bX_{n_1}\to \bX_{n_0}$, 
and let $n_2\geq n_1$
be another index, such that we have a multiplication 
$\on{mult}^{n_1}_{n_2,n_2}:\bX_{n_2}\times \bX_{n_2}\to \bX_{n_1}$, satisfying 
an obvious associativity with respect to $\on{mult}^{n_0}_{n_1,n_1}$.
For $m=1,2$ let us denote by $\on{act}^{n_m}_{j,i}$, $\on{act}^{n_m}_{k,j}$
the maps obtained by composing $\on{act}^{n_0}_{j,i}$ and 
$\on{act}^{n_0}_{k,j}$, respectively, with $\bX_m\to \bX_0$. 

We will construct a map $\bX_{n_2}\times \bV'_i\to \bV'_i$, which
amounts to a map $\bX_{n_2}^{\on{top}}\times \bV'_i\to \bV'_i$.
Let $v_j$ be an element in $\bW''\subset \bV_j$, and 
$h_{n_2}\in \bX_{n_2}^{\on{top}}$.
We claim that there exists an element, denoted 
$v'_j\in \bW''\subset \bV_j$, which is unique
modulo $\bW'$, satisfying
\begin{equation}  \label{cond}
\on{act}^{n_2}_{j,i}(h'_{n_2},v'_j)=
\on{act}^{n_1}_{j,i}(\on{mult}^{n_1}_{n_2,n_2}(h'_{n_2},h_{n_2}),v_j)\in
\bV_i,
\end{equation}
for any $h'_{n_2}\in \bX_{n_2}^{\on{top}}$. 

By assumption, every element $v_j\in \bW''$ can be written as
$\underset{a}\Sigma\, \on{act}^{n_2}_{k,j}(h_{n_2}^a,v_k^a)$ for
$h^a_{n_2}\in \bX_{n_2}^{\on{top}}$, $v_k^a\in \bV_k$.
For $h_{n_2}\in \bX_{n_2}^{\on{top}}$ as above we set
$$v'_j=\underset{a}\Sigma\,
\on{act}^{n_1}_{k,j}(\on{mult}^{n_1}_{n_2,n_2}(h_{n_2},h^a_{n_2}),v_k^a)\in 
\bW''\subset \bV_j.$$
It is easy to see that $v'_j$ satisfies \eqref{cond}.

\end{proof}

For $\BH$ as above we have a natural embedding $\on{triv}:\BVect\to 
\on{Rep}(\BH,\BVect)$,
corresponding to ``trivial'' representations.

\begin{cor}  \label{simple inv and coinv}
For $\BH$ satisfying ($*$), the functor $\on{triv}$ admits both
right and left adjoints.
\end{cor}

Note that in \propref{inv and coinv} a more general statement
is established.

\begin{proof}

First, from \secref{weakly strict} it follows the the functor
$\on{triv}:Vect\to \on{Rep}(\BH,Vect)$ admits right and left
adjoints, denoted $\Pi\mapsto \Pi^\BH$ and
$\Pi\mapsto \Pi_\BH$, respectively.

Therefore, using \propref{representations of pro-groups}, it
is enough to show that the functor
$\on{triv}:\on{Pro}(Vect)\to \on{Pro}(\on{Rep}(\BH,Vect))$
has left and right adjoints. But these are simply given by
sending $\Pi="\underset{\longleftarrow}{lim}"\, \Pi_i$ to
$\Pi_\BH\simeq "\underset{\longleftarrow}{lim}"\, (\Pi_i)_\BH$
and $\Pi^\BH\simeq "\underset{\longleftarrow}{lim}"\, (\Pi_i)^\BH$,
respectively.

\end{proof}

As every right adjoint, the functor $\Pi\mapsto \Pi^\BH$
is left-exact, and similarly, the functor 
$\Pi\mapsto \Pi_\BH$ is right-exact.

\begin{lem}  \label{exactness of Jacquet}
Assume that $\BH$ is the inverse limit of a weakly
surjective family of $\bH_i$, where each $\bH_i$ is a 
group-object in $\bSet$
isomorphic to a direct limit of $\bH_{i,j}$, with
each $\bH_{i,j}$ being a group-object of
$\on{Pro}(Set_0)$. Then the functor of coinvariants
$\on{Rep}(\BH,\BVect)\to \BVect$ is exact.
\end{lem}

\begin{proof}

According to \propref{representations of pro-groups} and 
\corref{simple inv and coinv},
each $\Pi\in \on{Rep}(\BH,\BVect)$ is an inverse limit of
$\Pi_k\in \on{Rep}(\BH,Vect)$, and
$\Pi_\BH\simeq "\underset{\longleftarrow}{lim}" (\Pi_k)_\BH$.
Therefore, it suffices to show that the functor of 
coinvariants is exact on $\on{Rep}(\BH,Vect)$. By the 
definition of the latter, we can replace $\BH$ by
one of its quotients $\bH_i$, which we will denote
by $\bH$.

However, the fact that functor $\Pi\mapsto \Pi_\bH$
is exact on the category $\on{Rep}(\bH,Vect)$ is well-known.
Indeed, if $\bH="\underset{\longrightarrow}{lim}" \bH_j$,
$\bH_j\in \on{Pro}(Set_0)$,
$$\Pi_\bH\simeq \underset{j}{\underset{\longrightarrow}{lim}}\, \Pi_{\bH_j},$$
but the functor $limInd$ is exact on $Vect$, and
the functor $\Pi\to \Pi_{\bH_j}$ is exact on
$\on{Rep}(\bH_j,Vect)$, since $\bH^{\on{top}}_j$ is a compact group.

\end{proof}

\ssec{}   \label{** condition}

Consider now the category $\BSet$ with its pseudo-action
on $\BVect$. The main object of study of this paper is
the category of representations $\on{Rep}(\BG,\BVect)$ of
a group-object $\BG\in \BSet$ in $\BVect$.
For brevity, we will denote the category by $\on{Rep}(\BG)$,
when no confusion is likely to occur.

\begin{lem} \label{pro-completeness}
The functor $limProj:\on{Pro}(\on{Rep}(\BG))\to \on{Rep}(\BG)$
is defined on the entire category and is exact.
\end{lem}

\begin{proof}

Recall (cf. \lemref{abelian}, with Ind replaced by Pro)
that the category $\BVect$ is closed under projective limits.
If $\Pi_i=(\BV_i,\rho_i)$ is an inverse system of objects
of $\on{Rep}(\BG)$, we define $\BV\in \BVect$ as 
$\underset{\longleftarrow}{lim}\, \BV_i$. It is easy to see
from the definitions that there exists an action
$\rho:\BG\times \BV\to \BV$, such that $(\BV,\rho)$
represents the projective limit 
$\underset{\longleftarrow}{lim}\, \Pi_i$. The exactness follows
from the fact that the functor
$limProj:\on{Pro}(\BVect)\to \BVect$ is exact.

\end{proof}

We will say that $\BH\in \BSet$ satisfies condition ($**$) if,
as an object of $\on{Ind}(\on{Pro}(\bSet))$, 
$\BH$ can be represented as $"\underset{\longrightarrow}{lim}"\, \BX_k$,
with $\BX_k\in \on{Pro}(\bSet)$ being weakly strict.

\medskip

As before, we have an obvious functor
$\on{triv}:\BVect\to \on{Rep}(\BH)$ corresponding to ``trivial''
representations.

\begin{prop} \label{inv and coinv}
The functor $\on{triv}:\BVect\to \on{Rep}(\BH)$ admits
a left adjoint, and when $\BH$ satisfies ($**$), also
a right adjoint.
\end{prop}

\begin{proof}

Let us first construct the left adjoint of $\on{triv}$.
Consider the covariant functor on the category $Vect$
that sends a vector space $\bV$ to
$\on{Hom}_{\on{Rep}(\BH)}(\Pi,\on{triv}(\bV))$. 
This functor is a subfunctor of
$\bV\mapsto \on{Hom}_{\BVect}(\Pi,\bV)$. Hence, by
\propref{prorep dom}, it is pro-representable.

Let us denote the resulting object of $\BVect$
by $\Pi_\BH$.
It is strightforward to check that for $\BV\in \BVect$, we
have a functorial isomorphism 
$\on{Hom}_{\on{Rep}(\BH)}(\Pi,\on{triv}(\BV))\simeq 
\on{Hom}_{\BVect}(\Pi_\BH,\BV)$.

\medskip

Now let us construct the right adjoint to $\on{triv}$.
Let us write $\BH\in \BSet$ as
$"\underset{\longrightarrow}{lim}"\, \BX_k$, where 
$\BX_k\in \on{Pro}(\bSet)$ are weakly strict. 

For a weakly strict object $\BX\in \on{Pro}(\bSet)$, 
$\BV,\BU\in \BVect$, and an 
action map $\phi:\BX\times \BV\to \BU$, consider
the kernel of $\phi$ as a functor on $\BVect$:
$$\on{ker}(\phi)(\BW)=\{\psi:\BW\to \BV\,|\, \phi\circ \psi=0\}.$$
We claim that this functor is representable. If this is so,
it is easy to see that the sought-for right adjoint of $\on{triv}$
is representable by 
$$(\BV,\rho)^\BH=\underset{k}{\underset{\longleftarrow}{lim}}\, 
\on{ker}\left(p-\on{act}:\BX_k\times \BV\to \BV\right),$$
where $\underset{\longleftarrow}{lim}$ is taken in the category 
$\BVect$, and $p$ is the obvious projection
map $\BX_k\times \BV\to \BV$.

To show the representability, we can assume that $\BU=\bU\in Vect$. Indeed, if
$\BU="\underset{\longleftarrow}{lim}"\, \bU_i$, then
$\on{Ker}(\phi)=\underset{i}{\underset{\longleftarrow}{lim}}\,
\on{ker}\left(\BX\times \BV\to \bU_i\right)$.
In the latter case, we can assume that 
$\BV="\underset{\longleftarrow}{lim}"\, \bV_j$,
and we have a compatible system of maps $\phi_j:\BX\times \bV_j\to \bU$.
By \secref{weakly strict}, $\on{ker}(\phi_j)\subset \bV_j$ is well-defined, and
it is easy to see that $"\underset{\longleftarrow}{lim}"\, 
\on{ker}(\phi_j)\in \BVect$
represents $\on{ker}(\phi)$.

\end{proof}

The main source of examples of such $\BG$, i.e., of group-objects
in $\BSet$, is provided by considering sets of 
points of algebraic groups with values in a two-dimensional local
field.

\ssec{}

Let $\bK$ be a local field, with the corresponding local ring
$\CO_\bK$. We will denote by $\pi$ a uniformizer of $\bK$.
Set $\bF=\bK((t))$, $\CO_\bF=\bK[[t]]$.

\medskip

Let $Sch^{ft}$ denote the category of separated schemes of finite type
over $\bK$. If $S$ is an object of $Sch^{ft}$, we will denote
by $S(\bK)$ the corresponding set of $\bK$-points. It is well-known
that $S(\bK)$ carries a natural locally compact totally disconnected
topology; therefore, as a topological space,
$S(\bK)\simeq \bS^{\on{top}}$
for a canonically defined locally compact object $\bS\in \bSet$. 

Hence, we obtain a functor
$S\mapsto \bS:Sch^{ft}\to \bSet$, and also the functors
$\on{Pro}(Sch^{ft})\to \on{Pro}(\bSet)$,
and $\on{Ind}(\on{Pro}(Sch^{ft}))\to \BSet$.

In particular, any affine scheme (not necessarily of finite type) 
over $\bK$ defines an object of $\on{Pro}(Sch^{ft})$, and hence,
an object of $\on{Pro}(\bSet)$. 
In addition, for any scheme of finite type $S$, the corresponding scheme 
of arcs $S[[t]]$ is naturally an object of $\on{Pro}(Sch^{ft})$:
$$S[[t]]\simeq "\underset{\longleftarrow}{lim}"\,S[t]/t^i.$$
We will denote the corresponding object of $\on{Pro}(\bSet)$ by
$\bS[[t]]$.

If $S$ is smooth, the maps in this family defining $S[[t]]$
are fibrations into affine spaces; therefore the corresponding maps 
$\bS[t]/t^j\to \bS[t]/t^i$ are weakly surjective. Hence, if $S$ is
smooth, the object $\bS[[t]]\in \on{Pro}(\bSet)$ is weakly strict.

\medskip

For a scheme $S'$ over $\bF$, we define its "restriction of
scalars" from $\bF$ to $\bK$ as a functor on the category 
of schemes over $\bK$ by 
$S\mapsto \on{Hom}_{\bF}(S\underset{\bK}\otimes \bF,S')$.
If $S'$ is of finite type and affine, then by
embedding it into an affine space one shows that the
above functor is ind-representable by an ind-scheme, which is
a direct limit of affine schemes under closed embeddings.
By taking $S'=S\underset{\bK}\otimes \bF$ for $S$ an affine
scheme of finite type over $\bK$, we obtain an object of 
$\on{Ind}(\on{Pro}(Sch^{ft}))$ that will be denoted by $S((t))$. 
The resulting object of $\BSet$ will be denoted by $\bS((t))$ or $\BS$.

By applying the functor of iterated inductive and projective limits
$\BSet\to Set$, we obtain from $\BS$ (resp., $\bS[[t]]$) the set,
which is tautologically identified with the set $S(\bF)$ of $\bF$-points
of $S$ (resp., $S(\CO_\bF)$--the set of $\CO_\bF$-points of $S$).

\ssec{}  \label{group notation}

If $G$ is a smooth linear algebraic group over $\bK$, by applying the
functor $G\mapsto \bG$ we obtain the corresponding
group-object in $\bSet$. In particular, we
can consider the category of representations
$\on{Rep}(\bG,Vect)$, which is tautologically
equivalent to the category of smooth 
representations of the locally compact group $G(\bK)$.

For a non-negative integer $i$, let us denote by
$G^i$ the congruence subgroup of $G[[t]]$, i.e.,
the kernel of $G[[t]]\to G[[t]]/t^i$; in particular,
$G^0=G[[t]]$. Let $\bG^i$ be the corresponding object of 
$\on{Pro}(\bSet)$.
A subgroup $\BH$ of $\bG[[t]]$ will be called {\it thick} if it 
contains $\bG^i$ for some $i$ and equals the preimage of 
a closed subgroup of $\bG[[t]]/\bG^i$ (we are slightly abusing
the terminology by identifying $\bG[[t]]/\bG^i$ with the
corresponding locally compact group).

For a thick  $\BH \subset \bG[[t]]$ we can consider the corresponding
categories $\on{Rep}(\BH,Vect)$ and $\on{Rep}(\BH,\BVect)$.
As was remarked above, $\bG[[t]]\in \on{Pro}(\bSet)$ is weakly
strict, and so are the groups $\bG^i$. From this it is easy
to see that any thick  subgroup $\BH\subset \bG[[t]]$ satisfies 
condition ($*$) of \secref{condition star}.

\medskip

Finally, for an algebraic group $G$ as above, we can
consider $\BG$ (sometimes also denoted $\bG((t))$), 
which is a group-object in
$\BSet$ and the corresponding category $\on{Rep}(\BG,\BVect)$,
which we will denote for brevity by $\on{Rep}(\BG)$. 

It is well-known that the ind-scheme $G((t))$ can be represented as
a direct limit under closed embeddings of subschemes, each of which
is stable under (both left and right) multiplication by $G[[t]]$, and is a
principal $G[[t]]$-bundle over a scheme of finite type. (In fact,
the above family of subschemes is obtained by taking the preimages
of finite-dimensional subschemes of the affine Grassmannian of $G$,
i.e., $\Gr_G=G((t))/G[[t]]$.) This
implies, in particular, that $\BG$ satisfies condition ($**$),
cf. \secref{** condition}.

\medskip

Let us denote by $\BVect^{G(\bF)}$ the category consisting of 
objects of $\BVect$ with an action of the abstract group $G(\bF)$.

\begin{lem}  \label{abstract group}
The natural forgetful functor $\on{Rep}(\BG)\to\BVect^{G(\bF)}$
is fully faithful.
\end{lem}

\begin{proof}

We have to show that if $(\BV_1,\rho_1)$ and $(\BV_2,\rho_2)$ are
two objects of $\on{Rep}(\BG)$, and $\BV_1\to \BV_2$ is a map
preserving the action of $G(\bF)$, then it is compatible with the
$\BG$-action.

This can be shown in the following general set-up: Let $\bV_1,\bV_2,\bW_1,\bW_2$
be vector spaces, and let $\BX$ be
an object of $\on{Pro}(\bSet)$ endowed with action maps $\BX\times \bV_k\to 
\bW_k$,
$k=1,2$. Let $\bV_1\to \bV_2$, $\bW_1\to \bW_2$ be maps, such that the
square
$$
\CD
\BX^{\on{top}}\times \bV_1 @>>> \BX^{\on{top}}\times \bV_2 \\
@VVV  @VVV \\
\bW_1 @>>> \bW_2
\endCD
$$
commutes, where $\BX^{\on{top}}$ is the topological space
obtained from the corresponding object of $\on{Pro}(Top^{Hlctd})$ by 
taking the projective limit.

Assume now that $\BX$ can be presented as $"\underset{\longleftarrow}{lim}"\,
\bX_i$, where the
maps $\bX_j\to \bX_i$ are such that the corresponding maps
$\bX_j^{\on{top}}\to \bX_i^{\on{top}}$ are surjective. Then it is easy to 
see that the square
$$
\CD
\BX\times \bV_1 @>>> \BX\times \bV_2 \\
@VVV  @VVV \\
\bW_1 @>>> \bW_2
\endCD
$$
commutes as well.

The above assumption is satisfied in our situation for $\BX$ being the 
object of $\on{Pro}(\bSet)$ corresponding to a subscheme of $G((t))$, 
obtained as a preimage of a closed  subscheme in $G((t))/G[[t]]$. The
required surjectivity follows from the fact that
the groups $G^i$ for $i>0$ are pro-unipotent.

\end{proof}

\ssec{Central extensions}   \label{cent ext}

Suppose now that $\wh{G}$ is a group-indscheme, which is a central
extension of $G((t))$ by the multiplicative group $G_m$, i.e.,
$$1\to G_m\to \wh{G}\to G((t))\to 1.$$  
In other words, $\wh{G}$ is a group-object in the category of
ind-schemes, such that if $G((t))="\underset{\longrightarrow}{lim}"\, X_k$,
and $X_k="\underset{\longleftarrow}{lim}"\, X_{k,l}$ with $X_{k,l}\in Sch^{ft}$,
then each $\wh{G}\underset{G((t))}\times X_k$ is a total space of
a $G_m$-torsor over $X_k$, and this torsor is pulled back from $X_{k,l}$
for some index $l$. In what follows we will assume that we have a 
splitting $G[[t]]\to \wh{G}$.

We will denote by $\widehat{\BG}$ the corresponding group-object
in $\BSet$, which is an extension of $\BG$ by $\bG_m$. 

Let $c$ be a character $G_m(\bK)\to \BC^*$. 
We will denote by $\on{Rep}_c(\wh{\BG})$ the category of
representations of $\widehat{\BG}$ with central character
$c$. In other words, the objects of this category are pairs
$\Pi=(\BV,\rho)$, where $\BV\in \BVect$, and $\rho$ is an
action map $\widehat{\BG}\times \BV\to \BV$, satisfying
the associativity and the unit axioms as above, and such that
the composite action
$$
\bG_m\times \BV\to \widehat{\BG}\times \BV\to \BV
$$
(where $\bG_m$ is viewed as an object of $\bSet\subset \BSet$)
corresponds to the above character.

\ssec{}

We propose the category $\on{Rep}(\BG)=\on{Rep}(\BG,\BVect)$ as a framework
for the study of representations of the group $G(\bF)$.
Let us explain why introducing pro-objects of $Vect$ appears
to be necessary.
For the remainder of this section, let us assume that 
$G$ is semi-simple, simply-connected and split.

\medskip

The first question to ask is whether the category 
$\on{Rep}(\BG)$ contains any objects $\Pi=(\bV,\rho)$, where
$\bV$ belongs to $Vect$. The answer is that such representations
are necessarily trivial (i.e., they lie in the
image of the functor $Vect\to\BVect\overset{triv}\to \on{Rep}(\BG)$),
for the same reason as 
why $\fp$-adic groups usually have no finite-dimensional
representations.

Indeed, suppose that $(\bV,\rho)$ is such a representation.
By \lemref{abstract group}, it is sufficient to prove that the corresponding
representation of the abstract group $G(\bF)$ on $\bV$ is trivial.

Consider the kernel $K$ of the action  $G(\bF)\times \bV\to \bV$. This
is a normal subgroup, and by definition, there exists an $i$ such that
$K\supset G^i(\bK)$. But then we claim that $K$ must coincide with
$G(\bF)$. Let $N$ be the maximal unipotent subgroup of $G$, and let
$N^i(\bK):=N(\bF)\cap G^i(\bK)$ be the corresponding congruence subgroup.
Then $N^i(\bK)\subset K$, but using the torus action and the normality
of $K$, we obtain that the entire $N(\bF)$ is contained in $K$. Again,
by normality, we obtain that all unipotent elements in $G(\bF)$ are
contained in $K$. However, it is known that for a split simply-connected
group, its set of field-valued points is generated by the subset of
unipotent elements.

\medskip

Another sense in which one may seek an alternative definition
of $G(\bF)$-representations is to consider the pseudo-action of
$\BSet$ on $\on{Ind}(Vect)=\on{Ind}(\on{Ind}(Vect_0))$.
We claim that (under the same assumption on $G$) all 
objects of $\on{Rep}(\BG,\on{Ind}(Vect))$ are again trivial.

\begin{proof}

As before, we have a fully faithful functor
$\on{Rep}(\BG,\on{Ind}(Vect))\to \on{Ind}(Vect)^{G(\bF)}$,
and it suffices to show that for any object 
$(\bV,\rho)$, $\bV\in\on{Ind}(Vect)$, the action of
the maximal unipotent group $N(\bF)$ on $\bV$ is trivial.
Obviously, we can replace $G$ by an $SL_2$ corresponding to 
some simple root; let $B\subset G$ be the corresponding
Borel subgroup, i.e., $N\simeq G_a$, and 
$B:=G_a\ltimes G_m$, where $G_m$ acts on $G_a$ by 
the square of the standard character.

Our $\bV$ is a direct limit $"\underset{\longrightarrow}{lim}"\bV_l$,
with $\bV_l\in Vect$. Fix an index $l$, and
it suffices to show that the action map $B(\bF)\times \bV_l\to \bV$
is trivial.

For a (not necessarily positive) integer $i$, 
let us denote by $N^i(\bK)$ the subgroup of $N(\bF)\simeq \bK((t))$
equal to $t^i\cdot \bK[[t]]$. If the action of $N(\bF)$ on $\bV_l$
is non-trivial, let $i$ be the minimal integer such that the
restriction of this action to $N^i(\bK)$ is trivial.
By assumption we have a non-trivial action map
$(N^{i-1}(\bK)/N^i(\bK)\simeq \bK)\times \bV_l\to \bV$.

Let now $j$ be a sufficiently large integer, so that 
the corresponding  congruence subgroup $(G_m)^j(\bK)$ acts trivially on $\bV_l$.
Take $i'=(i-1)-2j$ and consider now the action of $N^{i'}(\bK)$ on $\bV_l$.
Let $l'$ be a sufficiently large index such 
that the iteration of actions
$$N^{i'}(\bK)\times (G_m)^j(\bK)\times N^{i'}(\bK)\times (G_m)^j(\bK)
\times \bV_l\to \bV_{l'}$$
is well-defined. We will show that the action 
$N^{i-1}(\bK)/N^i(\bK)\times \bV_l\to \bV_{l'}$ is necessarily trivial,
which would be a contradiction. For that, it suffices to show that it  
is trivial on every element $v\in \bV_l$.
For every such $v$ there exists an integer $k$ such that the action of
$t^{i'}\cdot \pi^k\cdot \CO_\bK[[t]]\subset N^{i'}(\bK)$ on $v$ is trivial. 
Hence, for $g\in (G_m)^j(\bK)$ and $n\in t^{i'}\cdot \pi^k\cdot \CO_\bK[[t]]$
$$(n\cdot g\cdot n^{-1}\cdot g^{-1})\cdot v=v\in \bV_{l'}.$$
However, since $(G_m)^j(\bK)=1+t^j\cdot \bK[[t]]$,
the subset of $N(\bK)$ consisting of elements of the form
$(n\cdot g\cdot n^{-1}\cdot g^{-1})$, with $g$ and $n$ as above, equals
the entire $t^i\cdot\bK[[t]]$, in particular, it projects surjectively onto 
$N^{i-1}(\bK)/N^i(\bK)$. Therefore, for $n'\in N^{i-1}(\bK)$, we have
$n'\cdot v=v\in \bV_{l'}$, which is what we had to show.

\end{proof}

\section{The induction functor}   \label{induction}

\ssec{}

Let $G$ be a split reductive group over $\bK$, and let
$\BH$ be a thick subgroup of $\bG[[t]]$.
We have an obvious restriction functor
$r^{\BG}_{\BH}:\on{Rep}(\BG)\to \on{Rep}(\BH,\BVect)$.

Our goal in this section is define the functors
$\wt{i}^{\BG}_{\BH},i^{\BG}_{\BH}: \on{Rep}(\BH,\BVect)\to \on{Rep}(\BG)$,
such that $\wt{i}^{\BG}_{\BH}$ will be the right adjoint of 
$r^{\BG}_{\BH}$. 

We will have an injective functorial map 
$i^{\BG}_{\BH}(\Pi)\to \wt{i}^{\BG}_{\BH}(\Pi)$, and there is
a certain analogy between the functors $\wt{i}^{\BG}_{\BH}$ and 
$i^{\BG}_{\BH}$ and the functors of induction and compact induction 
in the theory of $\fp$-adic groups. When $\BH$ corresponds to a 
parahoric subgroup of $G[[t]]$, we will have an isomorphism 
$i^{\BG}_{\BH}\simeq \wt{i}^{\BG}_{\BH}$.

The construction of the functor $i^{\BG}_{\BH}$ makes sense
for any algebraic group $G$, but the construction of the functor 
$\wt{i}^{\BG}_{\BH}$ given below uses the fact that $G$ is reductive. 
However, we expect that the right adjoint to $r^{\BG}_{\BH}$
exists for any $G$.

\ssec{}  \label{induction functor}

To an object $\bX\in \on{Pro}(Set_0)$ we can attach the vector
space of locally constant $\BC$-valued functions,
denoted $\on{Funct}^{lc}(\bX)$.
Namely, if $\bX="\underset{\longleftarrow}{lim}"\, X_i$, 
$$\on{Funct}^{lc}(\bX)=\underset{\longrightarrow}{lim}\,\on{Funct}(X_i),$$
where the direct system is taken with respect to the 
pull-back maps between the spaces of functions. Of course,
$\on{Funct}^{lc}(\bX)$ identifies with the space of
locally constant functions on the topological space
$\bX^{\on{top}}$.

For any $\bX\in \bSet$ we define ${\mathbb Funct}^{lc}(\bX)\in \BVect$
by setting for
$"\underset{\longrightarrow}{lim}"\, \bX_j$, $\bX_j\in \on{Pro}(Set_0)$
$${\mathbb Funct}^{lc}(\bX)="\underset{\longleftarrow}{lim}"\, 
\on{Funct}^{lc}(\bX_j),$$
with respect to the restriction maps. We define the space 
$\on{Funct}^{lc}(\bX)\in Vect$ of locally constant functions on $\bX$
as $limProj({\mathbb Funct}^{lc}(\bX))$.

\medskip

If $\bX\in \bSet$ is locally compact (cf. \secref{set notation}),
we can introduce the vector space $\on{Funct}^{lc}_c(\bX)$, which can 
be called the space of locally constant functions with compact support.
One way to introduce it is as the space of locally-constant
compactly supported functions on $\bX^{\on{top}}$. Equivalently, if $\bX$ 
is represented as a direct limit as in \secref{set notation}, we have 
the natural ``extension by zero''
maps $\on{Funct}^{lc}(\bX_i)\to \on{Funct}^{lc}(\bX_j)$, and we set
$\on{Funct}^{lc}_c(\bX)=\underset{\longrightarrow}{lim}\, 
\on{Funct}^{lc}(\bX_i)$. 
Note that we always have an inclusion
$\on{Funct}^{lc}_c(\bX)\hookrightarrow \on{Funct}^{lc}(\bX)$.

Let $\bX$ be again locally compact, presented as a direct limit 
as in \secref{set notation}. If $\bX'\to \bX$ is map between objects 
of $\bSet$, we define the vector space $\on{Funct}^{lc}_{c,rel}(\bX')$ as
the inductive limit $\underset{\longrightarrow}{lim}\, 
\on{Funct}^{lc}(\bX'\underset{\bX}\times \bX_i)$.

\medskip

If $\bX\to \bY$
is a map in $\bSet$, we have the pull-back morphism 
$\on{Funct}^{lc}(\bY)\to \on{Funct}^{lc}(\bX)$, and
if this is a proper map between locally compact objects, we also have the 
morphism
$\on{Funct}^{lc}_c(\bY)\to \on{Funct}^{lc}_c(\bX)$.

Suppose now that $\bY^1,\bY^2\in \bSet$ are locally compact, and 
we have an action $\bX\times \bY^1\to \bY^2$ (in the sense of the canonical
tensor structure on $\bSet$), such that the map
$\bX\times \bY^1\to \bX\times \bY^2$ is proper. 
Then we obtain an action map 
$\bX\times \on{Funct}^{lc}_c(\bY^2)\to \on{Funct}^{lc}_c(\bY^1)$
(in the sense of the pseudo-action of $\bSet$ on $Vect$). 

For example, the above properness condition is always
satisfied if $\bX$ is a group-object acting on $\bY^1=\bY^2$. 

Note that the above action does not always extend onto
$\on{Funct}^{lc}(\bY)$.

\medskip

Let now $\bY$ be an object of $\on{Ind}(\bSet)$. 
We will say that $\bY$ is ``tame'' if it can be represented as 
$"\underset{\longrightarrow}{lim}"\, \bY_i$ 
such that $\bY_i\in \bSet$ are locally compact, and
the corresponding maps $\bY_i\to \bY_j$ are proper. 
If $\bY$ is ``tame'', we can attach to it the object 
$\on{Funct}^{lc}_c(\bY) \in \BVect$ as 
$\on{Funct}^{lc}_c(\bY)="\underset{\longleftarrow}{lim}"\,
\on{Funct}^{lc}_c(\bY_i)$, where the maps are again
given by restriction.

Suppose now that $\bY^1,\bY^2\in \on{Ind}(\bSet)$ are both
``tame'', $\bY^j_i="\underset{\longrightarrow}{lim}"\, \bY^j_i$ for $j=1,2$,
and let $\BX\times \bY^1\to \bY^2$
be an action of $\BX\in \BSet$ in the sense of
the pseudo-action of $\BSet$ on $\on{Ind}(\bSet)$. That is 
$\BX="\underset{\longrightarrow}{lim}"\, \BX_l$,
$\BX_l="\underset{\longleftarrow}{lim}"\, \bX_{l,k}$ and
the action is given by the maps
$\bX_{l,k}\times \bY_i^1\to \bY_{i'}^2$. 
We say that this action is proper if the above presentations
can be chosen so that the maps $\bX_{l,k}\times \bY_i^1\to 
\bX_{l,k}\times \bY_{i'}^2$ are proper. (This condition is
satisfied if $\BH$ is a group-object in $\BSet$ acting
on $\bY=\bY^1=\bY^2$.) 

If the action $\BX\times \bY^1\to \bY^2$ is proper we obtain an
action map $\BX\times \on{Funct}^{lc}_c(\bY_2)\to 
\on{Funct}^{lc}_c(\bY_1)$ in the sense of the canonical 
pseudo-action of $\BSet$ on $\BVect$.

\medskip

A little more generally, if $\bV$ is a vector space,
instead of complex-valued functions, we
can consider spaces of functions with values in $\bV$,
denoted $\on{Funct}^{lc}(\bX,\bV)$ and 
$\on{Funct}^{lc}_c(\bX,\bV)$, respectively.

\ssec{}

Let $i\geq 0$ be such that $\bG^i\subset \BH$. Consider the full subcategory of 
$\on{Rep}(\BH/\bG^i,Vect)\subset \on{Rep}(\BH,\BVect)$; 
we will first define the restrictions of 
the functors $i^\BG_\BH,\wt{i}^\BG_\BH$ to this subcategory.

Recall that there exists a strict ind-scheme of ind-finite-type
$G((t))/G^i$ ("strict" means that it can be presented 
as a direct limit of schemes with transition maps being
closed embedding).
Its existence, i.e., the ind-representability of
the corresponding functor, follows easily from the corresponding
fact for $\Gr_G=G((t))/G[[t]]$ (see, for example, the Appendix to \cite{Ga}).
As an object of $\on{Ind}(Sch^{ft})$ it 
carries an action of $G((t))\in \on{Ind}(\on{Pro}(Sch^{ft}))$
``on the left'' and a commuting action of 
$G([[t]]/t^i)\in Sch^{ft}$ ``on the right''.

Therefore, by applying the functor 
$S\mapsto \bS:Sch^{ft}\to \bSet$, we obtain 
a ``tame'' object, denoted $\BG/\bG^i$
in $\on{Ind}(\bSet)$, which carries the actions
of $\BG$ and $\bG[[t]]/\bG^i$.

For an object $\Pi=(\bV,\rho)\in \on{Rep}(\BH/\bG^i,Vect)$,
we obtain that $\on{Funct}^{lc}_c(\BG((t))/\bG^i,\bV)\in \BVect$
carries a natural $\BG$-action and a commuting
$\BH/\bG^i$-action.

The object of $\BVect$ underlying $i^\BG_\BH(\Pi)$ is set to
be 
\begin{equation}  \label{definition of i}
\left(\on{Funct}^{lc}_c(\BG/\bG^i,\bV)\otimes\mu(\BH/\bG^i)
\right)_{\BH/\bG^i},
\end{equation}
where $\mu(\BH/\bG^i)$ is the $1$-dimensional vector space
of left-invariant measures on the locally compact group 
$(\BH/\bG^i)^{\on{top}}$ (acted on naturally by $\BH/\bG^i$,
being a subspace of all measures on $(\BH/\bG^i)^{\on{top}}$).
The $\BG$-action on $\on{Funct}^{lc}_c(\BG/\bG^i,\bV)$
defines on $i^\BG_\BH(\Pi)$ a structure of an object of
$\on{Rep}(\BG)$.

This definition of $i^\BG_\BH(\Pi)$ can be rewritten as
follows. First, let us introduce the object 
$\BG/\BH\in \on{Ind}(\bSet)$.
Let us write $G((t))/G[[t]]$ as $"\underset{\longrightarrow}{lim}"\, S_k$,
$S_k\in Sch^{ft}$, and let $S_k^i$ be the preimage of $S_k$
in $G((t))/G^i$. Let $\bS_k$, $\bS_k^i$ be the corresponding
objects of $\bSet$. By construction $\bS_k^i$ carries an action 
of the groups $\BH\subset\bG[[t]]/\bG^i$, and we claim that the 
categorical quotient $\bS^\BH_k:=(\bS^i_k)/(\BH/\bG^i)\in \bSet$
is well-defined and is locally compact.
This follows for example from the fact that
$S_k^i\to S_k$ is a fibration locally trivial in the Zariski toplogy.
Let $\BG/\BH="\underset{\longrightarrow}{lim}"\, \bS_k^\BH$ be the
corresponding object of $\on{Ind}(\bSet)$; this object is
``tame'' and it evidently does not depend on the way we presented 
$G((t))/G[[t]]\in \on{Ind}(Sch^{ft})$ as a direct limit.

For $\Pi=(\bV,\rho)\in \on{Rep}(\BH/\bG^i,Vect)$,
let $\on{Funct}^{lc}_{c,rel}(\bS^i_k,\bV)$ be the space 
of locally-constant $\bV$-valued functions on $\bS^i_k$,
whose support is contained in the preimage of a 
compact subset of $\bS_k^\BH$ (see \secref{induction functor}).

We have:
\begin{equation} \label{integration}
\left(\on{Funct}^{lc}_c(\bS^i_k,\bV)\otimes\mu(\BH/\bG^i)
\right)_{\BH/\bG^i}\simeq
\left(\on{Funct}^{lc}_{c,rel}(\bS^i_k,\bV)\right)^{\BH/\bG^i},
\end{equation}
where the isomorphism is given by 
integration along the fibers of $\bS^i_k\to \bS^\BH_k$.

\medskip

The above isomorphism makes it clear that $i^\BG_\BH(\Pi)$,
as an object of $\on{Rep}(\BG)$, is independent of the choice of
the congruence subgroup $\bG^i$ contained in $\BH$.
In particular, we obtain a well-defined functor
$i^\BG_\BH:\on{Rep}(\BH,Vect)\to \on{Rep}(\BG)$. From \eqref{definition of i} 
we infer that  $i^\BG_\BH$ is right-exact, and from \eqref{integration} 
that it is also left-exact.

Using \propref{representations of pro-groups} we extend
the above functor $\on{Rep}(\BH,Vect)\to \on{Rep}(\BG)$ to a functor
$\on{Rep}(\BH,\BVect)\simeq \on{Pro}(\on{Rep}(\BH,Vect))\to
\on{Pro}(\on{Rep}(\BG))$, which is also exact. We extend it
further to an exact functor
$i^\BG_\BH:\on{Rep}(\BH,\BVect)\to \on{Rep}(\BG)$, using 
\lemref{pro-completeness}.

\medskip

It is easy to see that our functor $i^\BG_\BH$ is isomorphic
to the composition of two functors:
$i^\BG_{\bG[[t]]}:\on{Rep}(\bG[[t]],\BVect)\to \on{Rep}(\BG)$ and
$i^{\bG[[t]]}_\BH:\on{Rep}(\BH,\BVect)\to \on{Rep}(\bG[[t]],\BVect)$,
where the latter functor is defined by a similar induction procedure.

\ssec{}

Let us now define the functor
$\wt{i}^\BG_\BH:\on{Rep}(\BH,\BVect)\to \on{Rep}(\BG)$.
First, let us assume that $\BH$ is such that 
$\BG/\BH\in \on{Ind}(\bSet)$ is ind-compact, i.e.,
is a direct limit of compact objects of $\bSet$.
(E.g., this condition is verified for $\bG[[t]]$,
or more generally for any $\BH$ containing $\bI$, where
$I\subset G[[t]]$ is the Iwahori subgroup.
This follows from the fact that the affine flag
variety $G((t))/I$ is ind-proper, i.e., is a direct limit
of proper schemes of finite type.)

In this case we set $\wt{i}^\BG_\BH=i^\BG_\BH$.

\begin{prop} \label{adjointness}
If $\BG/\BH$ is ind-compact, the functor
$i^\BG_\BH$ is the right adjoint to
the restriction functor $r^\BG_\BH$.
\end{prop}

\begin{proof}

The proof mimics the proof of the usual
adjunction property for $\fp$-adic groups.

Let us first construct the adjunction map
$r^\BG_\BH\circ i^\BG_\BH\to \on{id}_{\on{Rep}(\BH,\BVect)}$.
By the definition of $i^\BG_\BH$, it is enough to 
construct a morphism
$r^\BG_\BH\circ i^\BG_\BH(\Pi)\to \Pi$ for an object
of $\on{Rep}(\BH/\bG^i,Vect)$ for some $i$.

For $\Pi=(\bV,\rho)\in \on{Rep}(\BH/\bG^i,Vect)$,
consider the canonical restriction map
$$\on{Funct}^{lc}_c(\BG/\bG^i,\bV)\otimes \mu(\BH/\bG^i)
\to \on{Funct}^{lc}_c(\BH/\bG^i,\bV)\otimes \mu(\BH/\bG^i),$$
which is bi-$\BH/\bG^i$-equivariant by construction.

Since $\on{Funct}^{lc}_c(\BH/\bG^i)\otimes \mu(\BH/\bG^i)$
identifies as a bi-module over $(\BH/\bG^i)^{\on{top}}$ with
the Hecke algebra of (compactly supported, locally constant)
measures on this group, we obtain a 
bi-$\BH/\bG^i$-equivariant map
$\on{Funct}^{lc}_c(\BH/\bG^i,\bV)\otimes \mu(\BH/\bG^i)\to \bV$.
By the $\BH/\bG^i$-equivariance on the right, we thus obtain a map
$\left(\on{Funct}^{lc}_c(\BG/\bG^i,\bV)\otimes 
\mu(\BH/\bG^i)\right)_{\BH/\bG^i}\to \bV$, as required.

\medskip

Let us now construct the second adjunction map
$(\BW,\rho')\to i^\BG_\BH\circ r^\BG_\BH(\BW,\rho')$ 
for $(\BW,\rho')\in \on{Rep}(\BG)$.
Using \propref{representations of pro-groups}, we can represent
$r^\BG_\BH(\BW)$ as an inverse limit of $\bW_i$, where each $\bW_i$
is a vector space underlying an object
of $\on{Rep}(\BH/\bG^{n_i},Vect)$ for some $n_i$.
Let $\BG/\BH="\underset{\longrightarrow}{lim}"\, \bS^\BH_k$ be as before,
and let $\bS_k^{n_i}$ be the preimage of $\bS^\BH_k$ in
$\BG/\bG^{n_i}$. 

Then the object of $\BVect$ underlying $i^\BG_\BH\circ r^\BG_\BH(\BW,\rho')$ 
is $$\underset{\underset{k,i}\longleftarrow}{"lim"}\, 
\left(\on{Funct}^{lc}_c(\bS_k^{n_i},\bW_i)\otimes \mu(\BH/\bG^{n_i})
\right)_{\BH/\bG^{n_i}}.$$

For every fixed $k$ and $i$, let an index $j$ be such that
the $\BG$-action on $\BW$ gives a map
$\on{act}_{j,i}:\bS_k^{n_j}\times \bW_j \to \bW_i$.
By further enlarging $j$, we may assume that this map
is compatible with the $\BH$-action.

We define a map 
$\bW_j\to \left(\on{Funct}^{lc}_c(\bS_k^{n_i},\bW_i)\otimes 
\mu(\BH/\bG^{n_i})\right)_{\BH/\bG^{n_i}}$ as follows.
First, the above action map gives rise to a map
$$\bW_j\to 
\left(\on{Funct}^{lc}(\bS_k^{n_j},\bW_i)\right)^{\BH/\bG^{n_j}}.$$

Now, from the fact that $\bS_k^\BH$ is compact and
isomorphism \eqref{integration}, we obtain
$$\left(\on{Funct}^{lc}(\bS_k^{n_j},\bW_i)\right)^{\BH/\bG^{n_j}}\simeq
\left(\on{Funct}^{lc}_c(\bS_k^{n_i},\bW_i)\otimes \mu(\BH/\bG^{n_i})
\right)_{\BH/\bG^{n_i}}.$$
By composing, we obtain the required morphism.

It is easy to check that the constructed map from
$\BW$ to the object of $\BVect$ underlying $i^\BG_\BH\circ r^\BG_\BH(\BW)$
respects the $\BG$-action. It is equally straightforward to see
that the two adjunction maps indeed give rise to the adjointness
of functors.

\end{proof}

Thus, to define the functor $\wt{i}^\BG_\BH$ in general, it suffices
to define the functor $\wt{i}^{\bG[[t]]}_\BH$, which is the right
adjoint to the restriction functor 
$r^{\bG[[t]]}_\BH: \on{Rep}(\bG[[t]],\BVect)\to \on{Rep}(\BH,\BVect)$.

Let $\bG^i$ be a congruence subgroup contained in $\BH$. We define the 
functor $$\wt{i}^{\bG[[t]]}_\BH: 
\on{Rep}(\BH/\bG^i,Vect)\to \on{Rep}(\bG[[t]]/\bG^i,Vect)$$ to equal
the corresponding functor defined for locally compact groups.

Explicitly, for $\Pi=(\bV,\rho)\in \on{Rep}(\BH/\bG^i,Vect)$,
$$\wt{i}^{\bG[[t]]}_\BH\simeq 
\left(\on{Funct}^{sm}(\bG[[t]]/\bG^i,\bV)\right)^{\BH/\bG^i},$$
where $\on{Funct}^{sm}(\bG[[t]]/\bG^i)$ is the space of 
functions on $\bG[[t]]/\bG^i$, smooth with respect to the action
of this group by left translations.

Note that for
$\Pi=(\bV,\rho)\in \on{Rep}(\BH/\bG^i,Vect)$ as above, the object of
$\BVect$ underlying $\wt{i}^\BG_\BH\circ \wt{i}^{\bG[[t]]}_\BH(\Pi)$ 
is $\underset{\underset{k}\longleftarrow}{lim}\, (\bW_k)^{\BH/\bG^i}$,
where each $\bW_k$ is a certain subspace of 
$\on{Funct}^{lc}(\bS^i_k,\bV)$, 
and $\bS^i_k$ is as in \eqref{definition of i}.
\medskip

The above functor $\on{Rep}(\BH/\bG^i,Vect)\to \on{Rep}(\bG[[t]]/\bG^i,Vect)$
extends to a functor 
$\wt{i}^{\bG[[t]]}_\BH:\on{Rep}(\BH,Vect)\to \on{Rep}(\bG[[t]],Vect)$.
Using \propref{representations of pro-groups}, from it
we obtain the  
functor $\wt{i}^{\bG[[t]]}_\BH:\on{Rep}(\BH,\BVect)\to \on{Rep}(\bG[[t]],\BVect)$,
which is the right adjoint to 
$r^{\bG[[t]]}_\BH: \on{Rep}(\bG[[t]],\BVect)\to \on{Rep}(\BH,\BVect)$;
and hence also the functor 
$\wt{i}^\BG_\BH:\on{Rep}(\BH,\BVect)\to \on{Rep}(\BG)$ with the 
desired adjointness
property.

\ssec{} \label{induction from Levi}

Consider the functor $r^\BG_\bG:\on{Rep}(\BG)\to \on{Rep}(\bG,\BVect)$
equal to the composition of 
$r^\BG_{\bG[[t]]}:\on{Rep}(\BG)\to \on{Rep}(\bG[[t]],\BVect)$
and the functor $\BV\mapsto \BV_{\bG^1}:
\on{Rep}(\bG[[t]],\BVect)\to \on{Rep}(\bG,\BVect)$. Note that
by \lemref{exactness of Jacquet}, $r^\BG_\bG$ is exact.
Its right adjoint, which we will denote by 
$i^\BG_\bG$ is the composition of the ``obvious''
functor $\on{Rep}(\bG,\BVect)\to \on{Rep}(\bG[[t]],\BVect)$ coming
from the homomorphism $\bG[[t]]\to \bG$ and the functor
$i^\BG_{\bG[[t]]}$ studied above.

\medskip

More generally, let $P\subset G$ be a parabolic,
with the Levi quotient $M$, and let
$I_P\subset \bG[[t]]$ be the corresponding
parahoric subgroup.
(For $P=B$ we will denote $I_B$ simply by $I$,
and $M$ by $T$).
Let $\bI_P$ (resp., $\bP$, $\bM$) be the corresponding
group-objects of $\on{Pro}(\bSet)$ (resp.,
$\bSet$).
 
In a similar fashion we obtain a pair of 
mutually adjoint functors 
$r^\BG_\bM:\on{Rep}(\BG)\to \on{Rep}(\bM,\BVect)$ and 
$i^\BG_\bM:\on{Rep}(\bM,\BVect)\to \on{Rep}(\BG)$.

\medskip

Let now $\wh{G}$ be a central extension of $G((t))$
by means of $G_m$ (cf. \secref{cent ext}), and let 
$\on{Rep}_c(\wh{\BG})$ be the corresponding category of 
representations.
Since we are given a splitting of $\wh{\BG}$ 
over $\bG[[t]]$, and hence, over $\BH$, we have an obvious restriction functor
$r^{\wh{\BG}}_{\BH}:\on{Rep}_c(\wh{\BG})\to\on{Rep}(\BH,\BVect)$.

By repeating the construction of the previous subsections,
we obtain the functors $i^{\wh{\BG}}_\BH:\on{Rep}(\BH,\BVect)\to 
\on{Rep}_c(\wh{\BG})$, and $\wt{i}^{\wh{\BG}}_\BH:\on{Rep}(\BH,\BVect)\to 
\on{Rep}_c(\wh{\BG})$, such that $\wt{i}^{\wh{\BG}}_\BH$ is the right adjoint
of $r^{\wh{\BG}}_\BH$, and $i^{\wh{\BG}}_\BH\simeq \wt{i}^{\wh{\BG}}_\BH$ when 
$\BH$ contains $\bI$.

We will denote by $r^{\wh{\BG}}_\bG$, 
$i^{\wh{\BG}}_\bG$ (resp., $r^{\wh{\BG}}_\bT$, $i^{\wh{\BG}}_\bT$) 
the corresponding functors between $\on{Rep}_c(\wh{\BG})$ and 
$\on{Rep}(\bG,\BVect)$
(resp., $\on{Rep}(\bT,\BVect)$).

\ssec{}
Next we will establish an analogue of Bernstein's geometric
lemma, which describes the composition of the functors
$i^{\wh{\BG}}_\bT$ and $r^{\wh{\BG}}_\bT$, cf. \cite{Be}.

Let $\Lambda$ be the lattice of co-weights of the maximal torus
$T$ of $G$, and $W$--the Weyl group. When we restrict $\wh{G}$
to $T((t))\subset G((t))$, the commutator defines a map
$$T((t))/T[[t]]\times T[[t]]\to G_m,$$
which factors through $T((t))/T[[t]]\twoheadrightarrow \Lambda$
and $T[[t]]\twoheadrightarrow T$. In other words, we have a map
$$\Lambda\times T\to G_m,$$ 
which defines a pairing $Q:\Lambda\otimes \Lambda\to \BZ$. 
This pairing is $W$-invariant,
since the central extension of $T((t))$ was induced from that of $G((t))$.
For $\lambda\in \Lambda$
we will denote by $\phi_{cQ}(\lambda)$ the character
$\bT\to \BC^*$ equal to
$$\bT\overset{Q(\lambda,\cdot)}\longrightarrow 
\bG_m\overset{c}\to \BC^*.$$

Recall now the affine flag
scheme of $G$, which is a strict ind-scheme, equal by definition 
to $G((t))/I$ (where $I$ is the Iwahori subgroup), and denoted $\on{Fl}_G$. 
Its existence, i.e., the ind-representability of the corresponding
functor, follows easily from the case of $\Gr_G$, cf. \cite{Ga}.

Recall also that the set of $I$-orbits on $\on{Fl}_G$ identifies
naturally with $W_{aff}\simeq \Lambda\ltimes W$--the extended affine Weyl 
group of $G$. For $w\in W_{aff}$, let us denote
by $\on{Fl}^w_G$ the corresponding orbit and by
$\ol{\on{Fl}}{}^w_G$ its closure. Note that $W_{aff}$ is
naturally partially ordered and
$\ol{\on{Fl}}{}^w_G=\underset{w'\leq w}\cup\, \on{Fl}^w_G$.

Let $\check T$ be the Langlands dual torus of $T$ (over
$\BC$), which identifies with the set of unramified
characters of $\bT$. 
For $w\in W_{aff}$ let us denote by 
$w(\rho_{aff})-\rho_{aff}$ the character of $T$
equal to the projection on $T$ of the sum of negative
affine roots which are turned positive by the action of
$w^{-1}$. Let $\mu_{w}$ denote the element in 
$\check T$ equal to the value of
$w(\rho_{aff})-\rho_{aff}:G_m\to \check T$ on $q\in \BC^*$,
where $q$ is the order of the residue field of our local
field $\bK$.

For $w\in W_{aff}$ and $\Pi\in \on{Rep}(\bT,Vect)$ we define a new 
representation $w\cdot \Pi$ by setting 
$$w\cdot \Pi:=\Pi^{\ol{w}}\otimes \phi_{cQ}(\lambda)\otimes \mu_{w},$$
where $w=\lambda\cdot \ol{w}$, $\lambda\in \Lambda$, $\ol{w}\in W$,
and $\Pi^w$ is obtained from $\Pi$ by twisting the $\bT$-action using 
$w$ viewed as an automorphism of $T$.

\begin{prop}  \label{geometric lemma}
For a representation $\Pi=(\bV,\rho)\in \on{Rep}(\bT,Vect)$,
the object $r^{\wh{\BG}}_\bT\circ i^{\wh{\BG}}_\bT(\Pi)$ can be canonically
written as 
$\underset{\underset{w\in W_{aff}}\longleftarrow} {"lim"}\, \bV_w$, 
$\bV_w\in \on{Rep}(\bT,Vect)$ in such a way that for
for $w'\leq w$ the map $\bV_w\to \bV_{w'}$ is a surjection, and the kernel
$\bV^w:=\on{ker}\left(\bV_w\to \underset{w'<w}\oplus\, \bV_{w'}\right)$ is 
isomorphic to $w\cdot \Pi$.
\end{prop}

\begin{proof}

Let $\bFl_G$ be the object of $\on{Ind}(\bSet)$ corresponding
to the ind-scheme $\on{Fl}_G$. Let us denote by $\bFl^w_G$ and
$\ol{\bFl}{}^w_G$ the corresponding objects of $\bSet$.

Let $\bI^0$ denote the kernel of the map $\bI\to \bT$.
Let $\bS^w$ (resp., $\ol{\bS}^w$) be the preimage of 
$\bFl^w_G$ (resp., $\ol{\bFl}{}^w_G$) in $\wh{\BG}/\bI^0$.
By construction, $r^{\wh{\BG}}_\bI\circ i^{\wh{\BG}}_\bT(\Pi)$
is the inverse limit of
$$\bW_w:=\left(\on{Funct}^{lc}_c(\ol{\bS}^w,\bV)\right)_{\bT\times \bG_m}
\simeq \left(\on{Funct}^{lc}(\ol{\bS}^w,\bV)\right)^{\bT\times \bG_m}.$$

Set $\bV_w:=\left(\bW_w\right)_{\bI^0}$. 
Since for $w'\leq w$, the restriction map
$\on{Funct}^{lc}_c(\ol{\bS}^w,\bV)\to
\on{Funct}^{lc}_c(\ol{\bS}^{w'},\bV)$ is surjective,
we obtain that $\bV_w\to \bV_{w'}$ are indeed surjective,
by the right-exactness of the functors $(\cdot )_{\bT\times \bG_m}$
and $(\cdot )_{\bI^0}$. 

Using \lemref{exactness of Jacquet}, 
we obtain that 
$$\bV^w:=\on{ker}(\bV_w\to \underset{w'<w}\oplus\, \bV_{w'})
\simeq \left(\left(\on{Funct}^{lc}_c(\bS^w,\bV)\right)_{\bT\times \bG_m}
\right)_{\bI^0}.$$

Let us choose a splitting $\bT\to \bB$, by means
of which $\bT$ becomes a subgroup of $\bI$; let
$g\in \on{Fl}^w_G(\bK)$ be the $\bT$-stable point, and let
$St(g)_{\bI}$ be the stabilizer of $g$ in $\bI$.
We obtain a homomorphism $St(g)_{\bI}\to \bT\times \bG_m$
and we have:
$$\left(\left(\on{Funct}^{lc}_c(\bS^w,\bV)\right)_{\bT\times \bG_m}
\right)_{\bI^0}\simeq
\left(i^{\bI}_{St(g)_{\bI}}\circ r^{St(g)_{\bI}}_{\bT\times \bG_m}(\Pi)
\right)_{\bI^0}.$$

\medskip

Observe that the character of $\bT$, 
corresponding to measures on the homogeneous space 
$\on{Fl}^w_G(\bK)\simeq \bI/St(g)_{\bI}$, equals $\mu_w$.

Write $w=\lambda\cdot \ol{w}$, $\lambda\in \Lambda$, $\ol{w}\in W$.
Observe now that the pull-back of $\Pi$ under the composition
$$\bT\to St(g)_{\bI}\to \bI\times \bG_m\to \bI\to\bT$$
is naturally isomorphic to $\Pi^{\ol{w}}$,
and the pull-back of the character $c:\bG_m\to \BC^*$ under
$\bT\to St(g)_{\bI}\to \bI\times \bG_m\to \bG_m$ is
$\phi_{cQ}(\lambda)$.

\medskip

This implies that $\bV^w \simeq w\cdot \Pi$ as $\bT$-representations.

\end{proof}

\ssec{}

One can formulate an analog of \propref{geometric lemma} 
describing the composition of the functors 
$r^{\wh{\BG}}_\bG\circ i^{\wh{\BG}}_\bG:\on{Rep}(\bG,\BVect)\to 
\on{Rep}(\bG,\BVect)$:

Set $\Gr_G:=G((t))/G[[t]]$, and recall that 
$G[[t]]$-orbits on $\Gr_G$ are in a natural bijection
with the partially ordered set $\Lambda^+$ of dominant
weights.

For every $\lambda\in\Lambda^+$, let $g\in \Gr_G^\lambda$
be a $T$-stable point, and let $St(g)_{G[[t]]}$ be its stabilizer
in $G[[t]]$, so that $\Gr_G^\lambda\simeq G[[t]]/St(g)_{G[[t]]}$.
Note that since $G^1$ is normal in $G[[t]]$, the quotient
$G^1\backslash  \Gr_G^\lambda$ is a $G$-homogeneous space 
isomorphic to $G/P^\lambda$ for a parabolic $P^\lambda\subset G$.
Let $M^\lambda$ be the Levi quotient of
$P^\lambda$, and let $\mu_\lambda$ be the character of $M^\lambda$
corresponding to measures on the homogeneous space $I_P/St(g)_{G[[t]]}(\bK)$.
Note also that for $\lambda$ as above, the character $\phi_{cQ}(\lambda)$
of $\bT$ is in fact well-defined as a character of $\bM^\lambda$.

Recall from the theory of $\fp$-adic groups that for a parabolic
$P$ with a Levi quotient $M$ we have a pair of mutually adjoint 
functors $r^\bG_\bM:\on{Rep}(\bG,Vect)\to \on{Rep}(\bM,Vect)$ and
$i^\bG_\bM:\on{Rep}(\bM,Vect)\to \on{Rep}(\bG,Vect)$.

\begin{prop}  \label{geometric lemma for Grassmann}
For a representation $\Pi=(\bV,\rho)\in \on{Rep}(\bG,Vect)$,
the object $r^{\wh{\BG}}_\bG\circ i^{\wh{\BG}}_\bG(\Pi)$ can be canonically
written as 
$\underset{\underset{\lambda\in \Lambda^+}\longleftarrow} {"lim"}\, \bV_\lambda$, 
$\bV_\lambda\in \on{Rep}(\bG,Vect)$ in such a way that for
$\lambda'\leq \lambda$ the map 
$\bV_\lambda\to \bV_{\lambda'}$ is a surjection, and the kernel
$\bV^\lambda:=\on{ker}(\bV_\lambda\to \underset{\lambda'<\lambda}\oplus\, 
\bV_{\lambda'})$ is 
canonically isomorphic to $i^\bG_{\bM^\lambda}\left(r^\bG_{\bM^\lambda}(\Pi)
\otimes \mu_\lambda\otimes \phi_{cQ}(\lambda)\right)$.
\end{prop}

The proof of this proposition is parallel to that of 
\propref{geometric lemma}.

\section{Examples}  \label{examples}

\ssec{}

Assume now that the group $G$ is split, simple and simply-connected.
In this case, a data of an extension $\wh{G}$ is equivalent 
to that of a $W$-invariant even symmetric bilinear form
$Q:\Lambda\otimes\Lambda\to \BZ$ (cf. Sect. 4 of \cite{BrDe}), and
we fix this form to be the minimal one, i.e., $\frac{1}{2\check h}Q_0$, 
where $Q_0$ corresponds to the Killing form, 
and $\check h$ is the dual Coxeter number. 

We have previously worked 
with a fixed character $\bG_m\to \BC^*$, but now we will consider 
all representations of the group $\wh{\BG}$. 
For a thick subgroup $\BH\subset \bG[[t]]$ we have the corresponding 
functors $i^{\wh{\BG}}_{\BH\times \bG_m}$, $r^{\wh{\BG}}_{\BH\times\bG_m}$ 
between $\on{Rep}(\BG)$ and $\on{Rep}(\BH\times \bG_m,\BVect)$.

Let $\Lambda_{aff}$ be the lattice $\Lambda\oplus\BZ$; which
identifies with the quotient of $\bT\times \bG_m$ by its maximal 
compact subgroup, and let $\BC[\Lambda_{aff}]$ be its 
group-algebra, viewed as a representation of $\bT\times \bG_m$.

Consider the object $\BV:=i^{\wh{\BG}}_{\bT\times \bG_m}(\BC[\Lambda_{aff}])\in
\on{Rep}(\wh{\BG})$,  studied by Kapranov in \cite{Kap}.
Let $\overset{\cdot\cdot}\sH_q$ be the modified Cherednik algebra of 
{\it loc.cit.} 2.3.3.
In \cite{Kap} it was shown that $\overset{\cdot\cdot}\sH_q$ injects into 
$\on{End}_{\on{Rep}(\wh{\BG})}(\BC[\Lambda_{aff}])$.

By combining the results of \cite{Kap} and \propref{adjointness}
we will establish the following result:

\begin{prop}  \label{Kapranov}
The map $\overset{\cdot\cdot}\sH_q\to \on{End}_{\on{Rep}(\wh{\BG})}(\BV)$ 
is an isomorphism.
\end{prop}

\begin{proof}

Let $\BV^{rat}$ be the object of $\on{Rep}(\wh{\BG})$ equal to
$i^{\wh{\BG}}_{\bT\times \bG_m}(\BC(\check T\times G_m))$, 
where $\BC(\check T\times G_m)$ is the field of rational functions on 
the torus $\check T\times G_m$, viewed as a $\bT\times \bG_m$-representation.
Note that by construction, both $\BV$ and $\BV^{rat}$
carry an action of the algebra $\BC[\Lambda_{aff}]$ by endomorphisms.

Using \propref{geometric lemma} and \propref{adjointness}
we obtain that
$$\on{Hom}_{\on{Rep}(\wh{\BG})}(\BV,\BV^{rat})\simeq 
\underset{\underset{w\in W_{aff}}\longrightarrow}{lim}\,
\on{Hom}_{\Lambda_{aff}-\on{mod}}(\bV_w,\BC(\check T\times G_m)),$$
where $\bV^w:=\on{ker}\left(\bV_w\to \underset{w'<w}\oplus\, \bV_{w'}\right)$
is isomorphic to $w\cdot \BC[\Lambda_{aff}]$. In particular, we
see that the restriction map
$\on{Hom}_{\on{Rep}(\wh{\BG})}(\BV^{rat},\BV^{rat})\to 
\on{Hom}_{\on{Rep}(\wh{\BG})}(\BV,\BV^{rat})$ is an isomorphism.

The subquotients $\on{Hom}_{\Lambda_{aff}-\on{mod}}(\bV^w,\BC(\check
T\times G_m))$ are all identified with $\BC(\check T\times G_m)$ as left
$\Lambda_{aff}$-modules, with the right $\Lambda_{aff}$-module
structure twisted by $w\cdot$ (see \propref{geometric lemma}).
Hence, we obtain a canonical direct sum decomposition
\begin{equation} \label{big Hom}
\on{End}_{\on{Rep}(\wh{\BG})}(\BV^{rat})\simeq
\on{Hom}_{\on{Rep}(\wh{\BG})}(\BV,\BV^{rat})\simeq \BC(\check T\times
G_m)\ltimes W_{aff}.
\end{equation}

\medskip

Therefore, using the main Theorem 3.3.8 of \cite{Kap}, it suffices to
check that the isomorphism \eqref{big Hom} coincides with the map
$$\overset{\cdot\cdot}\sH_q{}^{rat}\simeq 
\BC(\check T\times G_m)\ltimes W_{aff}\to 
\on{Hom}_{\on{Rep}(\wh{\BG})}(\BV,\BV^{rat})\to 
\on{End}_{\on{Rep}(\wh{\BG})}(\BV^{rat})
$$ of \cite{Kap}, Equation 3.3.7.

Since both isomorphisms preserve the ring structure, it suffices to
check that the generators of $\BC(\check T\times G_m)\ltimes W_{aff}$
over $\BC(\check T\times G_m)$, corresponding to the simple 
reflections under the two homomorphisms, act on $\BV^{rat}$ in the same way.

If $s$ is a simple reflection in $W_{aff}$, 
there exists a parahoric $I_s\subset G((t))$
such that the corresponding Levi quotient $M_s$ is a reductive group of
semi-simple rank $1$. As in \secref{induction from Levi} we have an induction
functor $i^{\wh{\BG}}_{\bM_s}:\on{Rep}(\bM_s,Vect)\to \on{Rep}(\BG)$, and
$\BV^{rat}\simeq i^{\wh{\BG}}_{\bM_s}\circ i^{\bM_s}_{\bT\times \bG_m}
(\BC(\check T\times G_m))$,
so that the endomorphism of $\BV^{rat}$ corresponding to $s$ via
both \eqref{big Hom} and the integral operator $\tau_s$ of \cite{Kap}
come from the corresponding endomorphisms of
$i^{\bM_s}_{\bT\times \bG_m}(\BC(\check T\times G_m))$. 

Therefore, we have reduced the question about the equality of two
endomorphisms of $\BV^{rat}$ to a similar question about endomorphisms
of $i^{\bM_s}_{\bT\times \bG_m}(\BC(\check T\times G_m))$ in the theory 
of $\fp$-adic groups. This reduces to the following (well-known)
calculation:

Let $G$ be a split reductive group of semi-simple rank $1$, and consider
the $G(\bK)$-representation $\bV:=i^\bG_\bT(\BC[\Lambda])$, which identifies
with the space of locally-constant compactly supported functions
on the quotient $G(\bK)/N(\bK)$, where $N$ is the maximal unipotent subgroup
of $G$. We can view $\bV$ as a $\check T$-family of principal series
representations, denoted $\bV_t, t\in \check T$. Let $\bV^{rat}$ be
the $\bG$-representation $i^\bG_\bT(\BC(\check T))$.
As above, we have
$$\BC(\check T)\ltimes W\simeq \on{End}_{G(\bK)}(\bV^{rat})\simeq 
\on{Hom}_{G(\bK)}(\bV,\bV^{rat}).$$
Consider the element $\tau_s$ of $\on{Hom}_{\bG}(\bV,\bV^{rat})$ corresponding
to the (unique) simple reflection in $W\subset \BC(\check T)\ltimes W$.
Then $\tau_s$ gives rise to a map $\bV_t\to \bV_{s\cdot t}$ defined
for $t$ belonging to an open subset of $\check T$, and the claim is that this
map is given by the meromorphic integral operator $f\mapsto f^{\tau_s}$ with
$$f^{\tau_s}(g)=\underset{n\in N(\bK)}\int\, f(g\cdot n\cdot s).$$

\end{proof}

\ssec{}

Let us go back to the situation, when the parameter $c$ is fixed and unramified. 
Let $\BV_{c}:=i^{\wh{\BG}}_\bT(\BC[\Lambda])\in \on{Rep}_c(\wh{\BG})$
and let $\overset{\cdot\cdot}\sH_{q,c}$ be the specialization of 
$\overset{\cdot\cdot}\sH_q$ at $c$.

\begin{cor}
We have an isomorphism 
$\overset{\cdot\cdot}\sH_{q,c}\simeq \on{End}(\BV_c)$.
\end{cor}

\begin{proof}

By applying \propref{adjointness} and \propref{geometric lemma} to 
$\BV$ and $\BV_c$, we obtain that
$\on{Hom}_{\on{Rep}(\wh{\BG})}(\BV,\BV)$ is a flat module over 
$\BC[\Lambda_{aff}]$, and hence over $\BC[\BZ]$, whose fiber
at $c\in \on{Spec}(\BC[\BZ])$ identifies with 
$\on{Hom}_{\on{Rep}_c(\BG)}(\BV_c,\BV_c)$.

In other words, $\on{Hom}_{\on{Rep}_c(\BG)}(\BV_c,\BV_c)$ is isomorphic
to the fiber of $\overset{\cdot\cdot}\sH_q$ at $c$, which is the same as 
$\overset{\cdot\cdot}\sH_{q,c}$.

\end{proof}

\ssec{}

The representation $\BV_c$ studied above is an analogue of principal series
representations. We will now introduce an object of $\on{Rep}(\BG)$ which
should be thought of as a cuspidal representation of $\BG$, although
at the moment we do not have a definition of cuspidality.

Let $\Pi$ be an irreducible cuspidal representation of $\bG$; and 
consider $i^{\wh{\BG}}_\bG(\Pi)\in \on{Rep}_c(\wh{\BG})$. 

\begin{lem}
$\on{End}_{\on{Rep}_c(\wh{\BG})}(i^{\wh{\BG}}_\bG(\Pi))\simeq \BC$.
\end{lem}

\begin{proof}

Using \propref{adjointness} we have
$\on{End}_{\on{Rep}_c(\wh{\BG})}(i^{\wh{\BG}}_\bG(\Pi))\simeq 
\on{Hom}_{\on{Rep}(\bG,\BVect)}(r^{\wh{\BG}}_\bG\circ i^{\wh{\BG}}_\bG(\Pi),\Pi)$.

We claim that the natural map 
$r^{\wh{\BG}}_\bG\circ i^{\wh{\BG}}_\bG(\Pi)\to \Pi$
is an isomorphism, which would imply the assertion of the lemma.
In fact, we claim that all the subquotients
$\bV^\lambda$ of \propref{geometric lemma for Grassmann} vanish 
except for $\lambda=0$. 

Indeed, by \propref{geometric lemma for Grassmann} each such
subquotient involves the functor $r^\bG_{\bM^\lambda}$ applied to
$\Pi$, which vanishes, since $\Pi$ was assumed to be cuspidal.

\end{proof}

\begin{conj}
The objects $i^{\wh{\BG}}_\bG(\Pi)$, for $\Pi$ being a cuspidal
representation of $\bG$, are irreducible.
\end{conj}

\ssec{}

Recall that an object $\Pi\in \on{Rep}(\bG,Vect)$ is called admissible
if for every open compact subgroup $\bH\subset \bG$, the vector
space $\Pi^\bH\simeq \Pi_\bH$ is finite-dimensional, i.e., 
belongs to $Vect_0$.

We can give an analogous definition in the case of
$\on{Rep}_c(\wh{\BG})$:

\begin{defn}
An object $\Pi\in \on{Rep}_c(\wh{\BG})$ is called admissible if for every 
thick subgroup $\BH\subset \bG[[t]]$ the object
$(\Pi)_\BH\in \BVect$ belongs, in fact, to $Vect$.
\end{defn}

It is easy to see that the principal series representations $\BV_c$ 
are not admissible. However, we have the following assertion: 

\begin{prop}
The representation $i^{\wh{\BG}}_\bG(\Pi)$, for $\Pi$ being a
cuspidal representation of $\bG$, is admissible.
\end{prop}

\begin{proof}

First, we can replace $\BH$ by a congruence subgroup $\bG^i$:
indeed, if $\BH\supset \bG^i$, then the statement
for $\bG^i$ would imply that for $\BH$.

Let $\bGr_G$ (resp., $\bGr^\lambda_G$, $\ol{\bGr}{}^\lambda_G$)
be the objects of $\on{Ind}(\bSet)$ and $\bSet$ corresponding
to $\Gr_G$, $\Gr^\lambda_G$ and $\ol{\Gr}{}^\lambda_G$, respectively.

Let $\bS^\lambda$ (resp., $\ol{\bS}{}^\lambda$) be the preimage
of $\bGr^\lambda_G$ (resp., $\ol{\bGr}{}^\lambda_G$) in $\wh{\BG}/\bG^1$.
As in \propref{geometric lemma}, the object
$r^{\wh{\BG}}_{\bG[[t]]}\circ i^{\wh{\BG}}_\bG(\Pi)\in \on{Rep}(\bG[[t]],\BVect)$ 
for $\Pi=(\bV,\rho)\in \on{Rep}(\bG,Vect)$ is the inverse limit 
over $\lambda\in \Lambda^+$ of
$\on{Funct}^{lc}_c\left(\ol{\bS}{}^\lambda,\bV\right)_{\bG\times \bG_m}$.

Set $\bW^\lambda:=
\on{Funct}^{lc}_c\left(\bS^\lambda,\bV\right)_{\bG\times \bG_m}$,
and we have to show that $\left(\bW^\lambda\right)_{\bG^i}\simeq 0$ for
all but finitely many $\lambda$'s.

\medskip

Let $g\in \Gr^\lambda_G(\bK)$ be a $\bT$-stable point and 
$\bSt(g)_{G[[t]]}$ its stabilizer in $\bG[[t]]$. 
By definition, we have a homomorphism
$\bSt(g)_{G[[t]]}\to \bG[[t]]\times \bG_m$ and
$$\bW^\lambda\simeq i^{\bG[[t]]}_{\bSt(g)_{G[[t]]}}\circ 
r^{\bSt(g)_{G[[t]]}}_{\bG\times \bG_m}(\Pi).$$
Therefore, as representations of $\bG[[t]]/\bG^i$,
$$\left(\bW^\lambda\right)_{\bG^i}\simeq 
i^{\bG[[t]]/\bG^i}_{\bSt(g)_{G[[t]]}/\bSt(g)_{G[[t]]}\cap \bG^i}
\left(\left(r^{\bSt(g)_{G[[t]]}}_{\bG\times\bG_m}
(\Pi)\right)_{\bSt(g)_{G[[t]]}\cap \bG^i}
\otimes \mu\right),$$
where $\mu$ is a character.

Note that as an object of $Vect$, 
$\left(r^{\bSt(g)_{G[[t]]}}_{\bG\times\bG_m}
(\Pi)\right)_{\bSt(g)_{G[[t]]}\cap \bG^i}$
is isomorphic to $(\Pi)_{\bH^i}$, where $\bH^i$ is the image of 
$\bSt(g)_{G[[t]]}\cap \bG^i$ under the homomorphism
$\bSt(g)_{G[[t]]}\to \bG[[t]]\to \bG$. Therefore, the assertion of the
proposition follows from the fact that for all but finitely many $\lambda$'s, 
the subgroup $\bH^i\subset \bG$ contains the unipotent radical of a non-trivial 
parabolic.

\end{proof}

\section{The Schwartz space on $\BG$}   \label{Schwartz space}

\ssec{}

Suppose now that $\bS$ is an object of $\bSet$, corresponding
to a smooth scheme of finite type $S$ over $\bK$. Then
it makes sense to consider the space of locally constant
compactly supported measures on $S(\bK)$, denoted $M(\bS)$.
To define it, we choose locally a top degree nowhere
vanishing differential form 
$\omega$ on $S$, which defines a measure $\mu(\omega)$ on $S(\bK)$. 
We say that a measure is locally constant if it can be obtained from
$\mu(\omega)$ by multiplication by a locally constant function.
One readily checks that this definition is independent of
the choice of $\omega$. 

Suppose now that $\phi:S_1\to S_2$ is a smooth map between 
smooth schemes. In this case, the operation of push-forward of 
constantly supported measures preserves the subspaces of 
locally constant ones, i.e., it defines a map
$\phi_!:M(\bS_1)\to M(\bS_2)$.

In particular, if an algebraic group $G$ acts on $S$, we obtain that
$M(\bS)$ is naturally on object of $\on{Rep}(\bG,Vect)$.

\medskip

For $S$ as above, consider now the object $\BS\in \BSet$. It appears
that there is no invariant way to assign to $\BS$ an object of
$\BVect$, which would be a replacement of locally constant 
compactly supported measures, and this is similar to the absence of a 
notion of D-module on $S((t))$, cf. \cite{AG}.

In this section we will study this phenomenon first when $S$ is
the affine space $A^n$, and then when $S$ is an affine algebraic 
group $G$.

\ssec{}

For any scheme $S$ which is isomorphic to a projective limit 
of smooth schemes of finite type $S_i$ with smooth transition
maps $S_i\to S_j$, we have $\bS\in \on{Pro}(\bSet)$, 
and we define $M(\bS)\in \BVect$ as
$M(\bS):="\underset{\longleftarrow}{lim}"\, M(\bS_i)$, where
the maps $M(\bS_j)\to M(\bS_i)$ for $j\geq i$ are the
push-forwards of measures. 

Recall that a lattice $L\subset \bK((t))^n$ is a finitely generated 
$\bK[[t]]$-submodule, which contains $t^i\cdot \bK[[t]]$ for some $i$.
The ``standard'' lattice is by definition $L_0=\bK[[t]]$.
By abuse of notation, we will denote by the same character $L$ the 
group-subscheme of $A^n((t))$ corresponding to a lattice $L$,
and by $\bL$ the corresponding object of $\on{Pro}(\bSet)$.
Since $L=\underset{i}{\underset{\longleftarrow}{lim}}\, (L/t^i\cdot L)$, we have 
a well-defined object $M(\bL)\in \BVect$.

\ssec{}

For a finite-dimensional vector space $H$ over $\bK$ let 
$\on{det}(H)=\Lambda^{\on{top}}(H)$ denote
its determant line. Let $\bH$ and $\on{det}(\bH)$ be the
corresponding objects of $\bSet$, and let $\mu(\on{det}(\bH))$ 
denote the $1$-dimensional $\BC$-vector space of Haar measures on 
$\on{det}(\bH)$. 
\footnote{Properly speaking, 
$\on{det}(H)$ is a super-vector space; however, in this paper it will 
appear only via $\mu(\on{det}(H))$, so the difficulties associated 
with the sign are irrelevant.} Of course, an element of $\mu(\on{det}(\bH))$
determines also a Haar measure on $\bH$.

Recall that for two lattices $L,L'\subset \bK((t))^n$ we can
assign their relative determinant line $\on{det}(L,L')$ so that 
$\on{det}(L,L'')\simeq
\on{det}(L,L')\otimes \on{det}(L',L'')$ and for $L\subset L'$,
$\on{det}(L,L')=\on{det}(L'/L)$, where the vector space 
$L'/L$ is, by definition, finite-dimensional. Let $\mu(\on{det}(\bL,\bL'))$ 
be the corresponding line of Haar measures.

\begin{lem}
For $L\subset L'$ we have a canonical morphism
$M(\bL')\to M(\bL)\otimes \mu(\on{det}(L,L'))$.
\end{lem}

\begin{proof}

Let $L''$ be a sublattice in $L$. By definition, for every such $L''$
we must construct a morphism $M(\bL'/\bL'')\to M(\bL/\bL'')\otimes
\mu(\on{det}(\bL,\bL'))$.
The required map is defined as a composition:
\begin{align*}
&M(\bL'/\bL'')\simeq \on{Funct}^{lc}_c(\bL'/\bL'')\otimes \mu(\on{det}(\bL'',\bL'))\to
\on{Funct}^{lc}_c(\bL/\bL'')\otimes \mu(\on{det}(\bL'',\bL'))\simeq \\
&\on{Funct}^{lc}_c(\bL/\bL'')
\otimes \mu(\on{det}(\bL'',\bL))\otimes
\mu(\on{det}(\bL,\bL'))\simeq M(\bL/\bL'')\otimes \mu(\on{det}(\bL,\bL')),
\end{align*}
where the arrow corresponds to the ordinary restriction of functions.

\end{proof}

Finally, we are ready to define the object $M(\BA^n)\in \BVect$, which
we propose as a candidate for the Schwartz space of functions on $A^n(\bF)$:
$$M(\BA^n):=\underset{L}{\underset{\longleftarrow}{lim}}\, 
M(\bL)\otimes \mu(\on{det}(\bL,\bL_0)),$$
where $\underset{\longleftarrow}{lim}$ is taken in the category 
$\BVect$, and the arrows are given by the lemma above.

\medskip

It is easy to see that the action of $A^n((t))$ on itself
by translations makes $M(\BA^n)$ an object of $\on{Rep}(\BA^n)$.

Recall also that the group-indscheme $GL_n((t))$, which acts naturally
on $A^n((t))$, has a canonical central extension $\wh{GL}_n$ by means
of $G_m$, whose $S$-points for a test-scheme $S$ are pairs 
$g\in \on{Hom}(S,GL_n((t)))$ and a trivialization 
of the line bundle $\on{det}(g\cdot L_0,L_0)$ on $S$.

\begin{thm}   \label{flat space}
The action of $GL_n((t))$ on $A^n((t))$ makes
$M(\BA^n)$ an object of $\on{Rep}(\wh{\BG L}_n)$, where
$\bG_m\subset \wh{\BG L}_n$ acts via the character
$\bG_m\to \BZ\overset{1\mapsto q}\longrightarrow \BC^*$.
\end{thm}

\begin{proof}

By construction, as an object of $\BVect$,
$$M(\BA^n)\simeq "\underset{\longleftarrow}{lim}"\,
M(\bL/\bL') \otimes \mu(\on{det}(\bL,\bL_0)),$$
where the inverse limit is taken over the partially ordered
set of pairs of lattices $L'\subset L$
with $(L'_1\subset L_1)\leq (L'_2\subset L_2)$ if and only if
$L_1\subset L_2$ and $L_1'\supset L'_2$.

For clarity, let us first define an action of the abstract group
$\wh{GL}_n(\bK)$ on
$M(\BA^n)$. For a pair of lattices $L'\subset L$ and
$g\in GL_n((t))(\bK)$, the action of $g$ defines an isomorphism
$M(\bL/\bL')\simeq M(g\cdot \bL/g\cdot \bL')$ and an isomorphism
$\on{det}(L,L_0)\simeq \on{det}(g\cdot L,g\cdot L_0)\simeq
\on{det}(g\cdot L,L_0)\otimes \on{det}(L_0,g\cdot L_0)$.
Hence, if we lift $g$ to an element of $\wh{GL}_n(\bK)$,
we obtain an isomorphism $\on{det}(L,L_0)\simeq \on{det}(g\cdot L,L_0)$,
i.e., we obtain a desired action.

\medskip

Let us now repeat this construction in order to obtain an action map
$\wh{\BG L}_n\times M(\BA^n)\to M(\BA^n)$.
Let us write $\wh{GL}_n$ as $"\underset{\longrightarrow}{lim}"\, S_k$,
and $S_k="\underset{\longleftarrow}{lim}"\, S_{k,l}$ with
$S_{k,l}\in Sch^{ft}$. Set $\BS_k$ (resp., $\bS_{k,l}$) to be
the corresponding objects of $\on{Pro}(\bSet)$ (resp., $\bSet$). 

\medskip

Recall that if $S$ is a scheme, there is a notion of an $S$-family
of lattices in $\bK((t))^n$, which is in fact the same as an $S$-point
of the affine Grassmannian of $GL_n$. If $L$ and $L'$ are two
$S$-families of lattices with $L'\subset L$, then the quotient
$L/L'$ is a vector bundle on $S$.

\medskip

For a pair of lattices $L'\subset L\subset \bK((t))^n$ and an index $k$,
using the action of $GL_n((t))$ on $\Gr_{GL_n}$, we obtain the
$S_k$-families of lattices that we will denote by $S_k\cdot L'\subset S_k\cdot L$.
Moreover, there exists another pair of lattices $L'_1\subset L_1$, thought of as
constant $S_k$-families, such that $L'_1\subset S_k\cdot L'$ and 
$S_k\cdot L\subset L_1$. Consider the quotients
$$H_{S_k}:=S_k\cdot L'/L'_1\subset  H'_{S_k}:=S_k\cdot L/L'_1\subset  
H''_{S_k}:=L_1/L'_1$$
as vector bundles on $S_k$. Note that both $H'_{S_k}/H_{S_k}$ and 
$H''_{S_k}$ are trivial
bundles with fibers $L/L'$ and $L_1/L'_1$, respectively. By the definition
of $\wh{GL}_n$, the line bundle $\on{det}(H''_{S_k}/H'_{S_k})$ is identified with
the trivial line bundle with fiber $\on{det}(L,L_1)$. Finally, 
there exists an index $l$, so that $H_{S_k}$ and $H'_{S_k}$, together 
with their embeddings into $H''_{S_k}$, come from vector bundles on $S_{k,l}$, 
which we will denote by $H_{S_{k,l}}$ and $H'_{S_{k,l}}$, respectively.

We need to construct an action map $\bS_{k,l}\times 
\left(M(\bL_1/\bL'_1)\otimes \mu(\on{det}(\bL_1,\bL))\right)\to 
M(\bL/\bL')$, which on the level of $\bK$-points amounts to the one 
constructed above. For that, by Zariski localizing $S_{l,k}$, we may assume 
that the vector bundle $H'_{S_{k,l}}$ on $S_{k,l}$ can be 
trivialized, i.e., $H'_{S_{k,l}}\simeq H'\times S_{k,l}$. 

Thus, we have a map $\bS_{k,l}\times \bH'\to \bL_1/\bL'_1$, such that
the corresponding map $\bS_{k,l}\times \bH'\to \bS_{k,l}\times \bL_1/\bL'_1$
is proper (cf. \secref{induction functor}). Therefore, it
defines an action map $$\bS_{k,l}\times \on{Funct}^{lc}_c(\bL_1/\bL'_1)\to 
\on{Funct}^{lc}_c(\bH').$$ By tensoring with $\mu(\on{det}(\bH'))$ we obtain
an action map $\bS_{k,l}\times (M(\bL_1/\bL'_1)\otimes 
\mu(\on{det}(\bL_1,\bL)))\to
M(\bH')$.

Similarly, we have a map $\bS_{k,l}\times \bH'\to \bL/\bL'$, and by integration
we obtain an action map $\bS_{k,l}\times M(\bH')\to M(\bL/\bL')$. The composition
\begin{align*}
&\on{Hom}(\bS_{k,l}\otimes M(\bH'),M(\bL/\bL'))\times 
\on{Hom}(\bS_{k,l}\otimes M(\bL_1/\bL'_1),M(\bH'))\to \\
&\on{Hom}((\bS_{k,l}\times \bS_{k,l})\otimes M(\bL_1/\bL'_1),M(\bL/\bL'))\to 
\on{Hom}(\bS_{k,l}\otimes M(\bL/\bL'),M(\bH))
\end{align*}
yields the desired action.

\end{proof}

\ssec{}

Let now $G$ be an algebraic group over $\bK$. Let $G((t))$ be
the corresponding loop group and $\wh{G}$ its central extension
$1\to G_m\to \wh{G}\to G((t))\to 1$, as in \secref{cent ext}.
Let us fix a character $c:\bG_m\to \BC^*$.
We will now define an object $M_c(\BG)\in \BVect$, which will
underly an object $(M_c(\BG),\rho)\in \on{Rep}_c(\wh{\BG})$.

For every integer $i$ consider the trivial representation
$\BC$ of the corresponding congruence subgroup $\bG^i$ and consider
$i^{\wh{\BG}}_{\bG^i}(\BC)\otimes \mu(\bG[[t]]/\bG^i)\in \on{Rep}_c(\wh{\BG})$,
where $\mu(\bG[[t]]/\bG^i)$ is the $1$-dimensional space of
left-invariant Haar measures on the group $(\bG[[t]]/\bG^i)^{\on{top}}$.

By the construction of the functor $i^{\wh{\BG}}_{\bG^i}$ via
compactly supported functions on $\wh{\BG}/\bG^i$, for $j\geq i$
we have the morphisms
$$i^{\wh{\BG}}_{\bG^j}(\BC)\otimes \mu(\bG[[t]]/\bG^j)\to 
i^{\wh{\BG}}_{\bG^i}(\BC)\otimes \mu(\bG[[t]]/\bG^i)$$ 
given by fiber-wise integration. 
Set $M_c(\BG):=\underset{\longleftarrow}{lim}\,\, i^{\BG}_{\bG^j}(\BC)$,
where $\underset{\longleftarrow}{lim}$ is taken in $\BVect$.

For example, it is easy to see that when $G\simeq A^n$ 
(and the central extension is trivial), the space $M(\BG)$ 
obtained in this way identifies canonically with $M(\BA^n)$ 
considered above.

\ssec{}

Since each $\bG^i$ is a normal subgroup in $\bG[[t]]$, 
the terms of the inverse system defining $M_c(\BG)$ 
carry a commuting
$\bG[[t]]$-action on the right, which is respected by the
arrows. Hence, $M_c(\BG)$ carries an additional  
$\bG[[t]]$-action ``on the right'', which commutes with the
action of $\wh{\BG}$ ``on the left''.

This $\wh{\BG}-\bG[[t]]$-module structure on $M_c(\BG)$ 
allows to reinterpret
the functor $i^{\wh{\BG}}_\BH$ introduced earlier:

Recall that the functor of tensor product
$Vect\times Vect\to Vect$ extends naturally to
$\BVect$: 
$$("\underset{\longleftarrow}{lim}"\, \bV_i)\otimes
("\underset{\longleftarrow}{lim}"\, \bW_j):=
"\underset{\longleftarrow}{lim}"\, (\bV_i\otimes \bW_j).$$

Let $\BH\subset \bG[[t]]$ be a thick subgroup,
and let $\Pi=(\BV,\rho)$ be an object of $\on{Rep}(\BH,\BVect)$.
Consider the tensor product $M_c(\BG)\otimes \BV\in \BVect$.
The diagonal action of $\BH$ makes it into an object of
$\on{Rep}(\BH,\BVect)$, which carries a commuting $\wh{\BG}$-action.
Hence, $\left(M_c(\BG)\otimes \BV\right)_\BH$ is naturally an object
of $\on{Rep}(\BG)$.

The following is straightforward from the definitions:

\begin{lem}
We have a natural isomorphism in $\on{Rep}(\BG)$:
$i^{\wh{\BG}}_\BH(\Pi)\simeq \left(M_c(\BG)\otimes \BV\right)_\BH$.
\end{lem}

\ssec{}

Since we think of $M_c(\BG)$ as the space of functions 
on the group $\BG$, it is natural to expect that the
$\bG[[t]]$-action on $M_c(\BG)$ considered above extends
to an action ``on the right'' of the entire group $\BG=\bG((t))$, 
corresponding to right translations.
The existence of such an action is given by the theorem below.

Let $\fg$ be the Lie algebra of $G$, let $GL_{\fg}((t))$
be the corresponding loop group, and let $\wh{GL}_{\fg}$
be its canonical central extension as in \thmref{flat space}.
Let $\wh{G}_0$ be the central extension of $G((t))$ induced
from $\wh{GL}_{\fg}$ by means of the adjoint action.
For example, if $G$ is simple and simply-connected, the extension
$\wh{G}_0$ corresponds to the pairing 
$\Lambda\otimes \Lambda\to \BZ$ given by the Killing form.

Let $\wh{G}{}'$ be the central extension of $G((t))$ equal to
the Baer sum of $\wh{G}_0$ and the original
extension $\wh{G}$. Let $c'$ be the character of 
$\bG_m\subset \wh{\BG}{}'$ equal to the inverse of 
the product of $c$ and $c_0$, where $c_0$ is the character
$\bG_m\to \BZ\overset{1\mapsto q}\longrightarrow\BC^*$.

\begin{thm}  \label{main}
We have a canonical action of $\wh{\BG}{}'$ on 
$M_c(\BG)$, with $\bG_m\subset \wh{\BG}{}'$ acting
by the character $c'$. This action extends the natural
action of $\bG[[t]]$ on $M_c(\BG)$ ``on the right'' 
and commutes with the action of $\wh{\BG}$ ``on the left''.
\end{thm}

\ssec{Proof of \thmref{main}}

Let us first construct an action of the abstract group
$\wh{G}{}'(\bK)$ on $M_c(\BG)$.
For an integer $i$ and a point $g\in G((t))(\bK)$, there exists
an integer $j$ such that $g^{-1}(\bG^j)g$ is contained in
$\bG^i$; therefore, the right multiplication map
$\bG((t))\times g\to \bG((t))$ descends to a well-defined map
$\bG((t))/\bG^j\times g\to \bG((t))/\bG^i$.

In particular, if we lift $g$ to an element of 
$\wh{G}(\bK)$, we obtain a map 
$$i^{\wh{\BG}}_{\bG^j}(\BC)\otimes \mu(\bG^i/g^{-1}(\bG^j)g)\to
i^{\wh{\BG}}_{\bG^i}(\BC),$$
which commutes with the left $\wh{\BG}$-action. 

\medskip

We claim now that a lift of $g$ to an element of $\wh{G}_0(\bK)$ define
an identification of the line $\mu(\bG^i/g^{-1}(\bG^j)g)$
with $\mu(\bG^i/\bG^j)$. Indeed, let $\fg^i$ be the Lie subalgebra in $\fg((t))$ 
corresponding to the congruence subgroup $\bG^i$, then
$$\mu(\bG^i/\bG^j)\simeq \mu(\on{det}(\fg^i/\fg^j));\,\,
\mu(\bG^i/g^{-1}(\bG^j)g)\simeq \mu(\on{det}(\fg^i/g^{-1}\cdot\fg^j)),$$
and
$$\on{det}(\fg^i/\fg^j)\simeq \on{det}(\fg^i/g^{-1}\cdot \fg^j)\otimes
\on{det}(g\cdot \fg^j,\fg^j)\simeq 
\on{det}(\fg^i/g^{-1}\cdot \fg^j)\otimes \on{det}(g\cdot \fg^0,\fg^0).$$

\medskip

Therefore, if we take the Baer product of the extensions
$\wh{\BG}$ and $\wh{\BG}_0$ we obtain an action of 
the group of $\bK$-points of $\wh{\BG}{}'$ on $M_c(\BG)$.
Therefore, 
by passing to inverses, we obtain on $M_c(\BG)$ an action 
of $\wh{G}{}'(\bK)$, commuting with the left action of
$\wh{\BG}$ and the prescribed value of the central character.

\medskip

The fact the constructed point-wise action gives rise to 
a well-defined action map
$\wh{\BG}{}'\times M_c(\BG)\to M_c(\BG)$ follows
by considering families, as in the proof of \thmref{flat space}.

Namely, if $\wh{G}\underset{G((t))}\times \wh{G}_0="\underset{\longrightarrow}{lim}"\, S_k$,
$S_k="\underset{\longleftarrow}{lim}"\, S_{k,l}$ with
$S_{k,l}\in Sch^{ft}$, for every pair of indices $i,k$ there 
exists a large enough index $j$, such that the
group-subscheme $\on{Ad}_{S_k}(G^j)\subset G((t))\times S_k$
is contained in $G^i\times S_k$. Moreover, the relative determinant 
line $\on{det}(\fg^i/\on{Ad}_{S_k}(\fg^j))$ is identified 
with the constant line bundle with fiber $\on{det}(\fg^j,\fg^i)$. 

We have the map 
$$(\wh{G}/G^j)\times S_k\simeq (\wh{G}\times S_k)/\on{Ad}_{S_k}(G^j)
\to \wh{G}/G^i,$$
which comes from a map $\wh{G}/G^j\times S_{k,l}\to \wh{G}/G^i$
defined for a sufficiently large index $l$. The resulting map
$\wh{G}/G^j\times S_{k,l}\to \wh{G}/G^i\times S_{k,l}$
is smooth over every finite-dimensional subscheme of $\wh{G}/G^i$. 

Integration along the fiber defines the desired map
$$\bS_{k,l}\times \left(i^{\wh{G}}_{\bG^j}(\BC)\otimes \mu(\on{det}(\fg^j,\fg^i))\right)\to
i^{\wh{G}}_{\bG^i}(\BC).$$

\end{document}